\documentclass[11pt]{article}

\usepackage{amssymb,diagrams}
\usepackage{epsfig,amsmath,latexsym}
\usepackage{amsfonts}

\newtheorem{theorem}{Theorem}[section]

\newtheorem{defi}{Definition}[section]
\newtheorem{lemma}[theorem]{Lemma}

\def\slfrac#1#2{\hbox{\kern.1em %
 \raise.5ex\hbox{\the\scriptfont0 #1}\kern-.11em %
 /\kern-.15em\lower.25ex\hbox{\the\scriptfont0 #2}}}

\newcommand{\pf}{\noindent{\bf Proof.~}}
\newcommand{\beq}{\begin{eqnarray}}
\newcommand{\eeq}{\end{eqnarray}}
\newcommand{\beql}[1]{\begin{eqnarray}\label{#1}}
\newcommand{\beqs}{\begin{eqnarray*}}
\newcommand{\eeqs}{\end{eqnarray*}}
\newcommand{\eqn}[1]{(\ref{#1})}

\newcommand{\RR}{{\mathbb R}}
\newcommand{\ZZ}{{\mathbb Z}}
\newcommand{\QQ}{{\mathbb Q}}

\newcommand{\DD}{{\mathbb D}}

\newcommand{\rr}{{\mathbb R}}
\newcommand{\zz}{{\mathbb Z}}
\newcommand{\qq}{{\mathbb Q}}

\newcommand{\dd}{{\mathbb D}}

\newcommand{\fd}{{\mathfrak d}}
\newcommand{\fg}{{\mathfrak g}}

\newcommand{\fj}{{\mathfrak j}}

\newcommand{\fo}{{\mathfrak o}}

\newcommand{\ft}{{\mathfrak t}}

\newcommand{\bb}{{\mathbf b}}

\newcommand{\be}{{\mathbf e}}

\newcommand{\bh}{{\mathbf h}}
\newcommand{\bm}{{\mathbf m}}
\newcommand{\bo}{{\mathbf 1}}

\newcommand{\bw}{{\mathbf w}}
\newcommand{\bx}{{\mathbf x}}
\newcommand{\by}{{\mathbf y}}
\newcommand{\bz}{{\mathbf z}}

\newcommand{\bA}{{\mathbf A}}
  
\newcommand{\bC}{{\mathbf C}}
\newcommand{\bD}{{\mathbf D}}
\newcommand{\bI}{{\mathbf I}}
\newcommand{\bJ}{{\mathbf J}}
\newcommand{\bK}{{\mathbf K}}
\newcommand{\bM}{{\mathbf M}}
\newcommand{\bN}{{N}}
\newcommand{\bO}{{\mathbf O}}
\newcommand{\bP}{{\mathbf P}}
\newcommand{\bQ}{{\mathbf Q}}
\newcommand{\bS}{{\mathbf S}}
\newcommand{\bU}{{\mathbf U}}
\newcommand{\bV}{{\mathbf V}}
\newcommand{\bW}{{\mathbf W}}
\newcommand{\bX}{{\mathbf X}}
\newcommand{\bZ}{{\mathbf Z}}

\newcommand{\sA}{{\mathcal A}}

\newcommand{\sD}{{\mathcal D}}
\newcommand{\sE}{{\mathcal E}}

\newcommand{\sL}{{\mathcal L}}
\newcommand{\sM}{{\mathcal M}}

\newcommand{\sR}{{\mathcal R}}

\newcommand{\bsq}{\vrule height .9ex width .8ex depth -.1ex}

\makeatletter
\def\@sect#1#2#3#4#5#6[#7]#8{\ifnum #2>\c@secnumdepth
     \def\@svsec{}\else
     \refstepcounter{#1}\edef\@svsec{\csname the#1\endcsname.\hskip .75em }\fi
     \@tempskipa #5\relax
      \ifdim \@tempskipa>\z@
        \begingroup #6\relax
          \@hangfrom{\hskip #3\relax\@svsec}{\interlinepenalty \@M #8\par}%
        \endgroup
       \csname #1mark\endcsname{#7}\addcontentsline
         {toc}{#1}{\ifnum #2>\c@secnumdepth \else
                      \protect\numberline{\csname the#1\endcsname}\fi
                    #7}\else
        \def\@svsechd{#6\hskip #3\@svsec #8\csname #1mark\endcsname
                      {#7}\addcontentsline
                           {toc}{#1}{\ifnum #2>\c@secnumdepth \else
                             \protect\numberline{\csname the#1\endcsname}\fi
                       #7}}\fi
     \@xsect{#5}}
\def\@begintheorem#1#2{\it \trivlist \item[\hskip \labelsep{\bf #1\ #2.}]}

\def\plain{plain}\ifx\fmtname\plain\csname fi\endcsname
     
     \input docstrip
     \preamble

     Do not distribute the stripped version of this file.
     The checksum in the header refers to the documented version.

     \endpreamble
     \generateFile{here.sty}{t}{\from{here.doc}{}}
     \endinput
\fi
\ifcat a\noexpand @\let\next\relax\else\def\next{%
    \documentstyle[here,doc]{article}\MakePercentIgnore}\fi\next
\ifx\@Hxfloat\@Hundef\else\expandafter\endinput\fi
\let\@Hxfloat\@xfloat
\def\@xfloat#1[{\@ifnextchar{H}{\@HHfloat{#1}[}{\@Hxfloat{#1}[}}
\def\@HHfloat#1[H]{%
\expandafter\let\csname end#1\endcsname\end@Hfloat
\vskip\intextsep\vbox\bgroup\def\@captype{#1}\parindent\z@
\ignorespaces}
\def\end@Hfloat{\egroup\vskip \intextsep}

\makeatother

\catcode`\@=11

\catcode`@=12

\begin{document}


\begin{center}
{\Large {\bf Apollonian Circle Packings: Geometry and Group Theory\\
 III. Higher Dimensions}} \\

\vspace{1.5\baselineskip}
{\em Ronald L. Graham} \\
Department of  Computer Science and Engineering \\
University of California at San Diego,
La Jolla, CA 92093-0114 \\
\vspace*{1.5\baselineskip}

{\em Jeffrey C. Lagarias} \\
Department of Mathematics \\
University of Michigan, 
Ann Arbor, MI 48109--1109 \\
\vspace*{1.5\baselineskip}

{\em Colin L. Mallows} \\
Avaya Labs, Basking Ridge, NJ 07920 \\
\vspace*{1.5\baselineskip}

{\em Allan R. Wilks} \\
AT\&T Labs, Florham Park, NJ 07932-0971 \\
\vspace*{1.5\baselineskip}

{\em Catherine H. Yan}
\footnote{ 
Partially supported by NSF grants DMS-0070574, DMS-0245526 and a Sloan
Fellowship. This author is also affiliated with Dalian University of
Technology, China.}\\
Department of Mathematics \\
Texas A\&M University, 
College Station, TX 77843-3368\\
\vspace*{1.5\baselineskip}
(March 10, 2005  version) \\
\vspace*{1.5\baselineskip}
{\bf ABSTRACT}
\end{center}
This paper gives $n$-dimensional analogues of the
Apollonian circle packings in
parts I and II. Those papers  considered circle
packings described in terms of their Descartes
configurations, which are sets of four mutually
touching circles. They studied  packings 
that had integrality properties in
terms of the curvatures and centers of
the circles. Here we consider collections
of $n$-dimensional Descartes configurations,
which consist of $n+2$ mutually touching spheres.

We work in the 
space $\sM_{\dd}^n$ 
of all $n$-dimensional oriented Descartes configurations
parametrized 
in a coordinate system, ACC-coordinates, as those $(n+2) \times (n+2)$
real matrices $\bW$ with $\bW^T \bQ_{D,n} \bW = \bQ_{W,n}$ where 
$Q_{D,n} = x_1^2 +... + x_{n+2}^2 - \frac{1}{n}(x_1 +\cdots + x_{n+2})^2$ 
is the $n$-dimensional Descartes quadratic form, 
$Q_{W,n} = -8x_1x_2 + 2x_3^2 + \cdots + 2x_{n+2}^2$,
and $\bQ_{D,n}$ and $\bQ_{W,n}$
are their corresponding symmetric matrices.
 On the parameter space $\sM_{\dd}^n$
of ACC-matrices the  group
$Aut(Q_{D,n})$ acts on the left, and $Aut(Q_{W,n})$ acts on the right.
Both these  groups are  isomorphic to the 
$(n+2)$-dimensional Lorentz group $O(n+1,1)$, and give two
different ``geometric'' actions. The right action of $Aut(Q_{W,n})$ 
(essentially) corresponds to M\"{o}bius transformations acting 
on the underlying
Euclidean space $\rr^n$ while the left action of $Aut(Q_{D,n})$ is defined
only on the parameter space $\sM_{\dd}^n$.

We introduce  $n$-dimensional analogues of the 
Apollonian group, the dual Apollonian group and the
super-Apollonian group. These are finitely generated
groups in  $Aut(Q_{D,n})$, with 
the following integrality properties: 
the dual Apollonian group consists of
integral matrices in all dimensions, while the other two 
consist of rational matrices, with denominators having prime
divisors drawn from a finite set $S$ depending on the dimension.
We show that the  Apollonian group and the
dual Apollonian group are finitely presented, and are Coxeter groups.
We define an {\em Apollonian cluster ensemble} 
to be any orbit under the Apollonian group,
with similar notions for
the other two groups. We determine in which dimensions there exist 
rational Apollonian cluster ensembles
(all curvatures rational) and 
strongly rational Apollonian sphere ensembles (all ACC-coordinates rational).
\vspace*{1.5\baselineskip}

\setcounter{page}{1}

%
%
%
\setlength{\baselineskip}{1.0\baselineskip}

\section{Introduction}
\setcounter{equation}{0}
In part I we considered  Apollonian circle packings, and
observed that  there exist such packings where the circles
all had integer curvatures and the quantities
(curvature)$\times$(circle center) had integer 
entries. We gave an explanation for this phenomenon,
in terms of the  Descartes configurations
in the packing. A {\em Descartes configuration} is
a configuration of four mutually tangent circles,
with all six tangency points distinct.
We introduced  the space $\sM_{\dd}$ of all ordered, oriented Descartes
configurations, parametrized in a 
coordinate system expressed in terms of the
curvatures and centers of the circles in the
packing. 

The explanation of the integrality properties was as follows.
This configuration  space of ordered, oriented Descartes configurations
was shown to be 
a principal homogeneous space for the  Lorentz  group $O(3, 1)$,
a six-dimensional real Lie group which we 
identified with the real automorphism group 
$Aut(Q_{D})$ of a quaternary quadratic form
$$
Q_{D}(x_1, x_2, x_3, x_4) = x_1^2 + x_2^2 + x_3^2 + x_4^2 -
\frac{1}{2} (x_1 + x_2 + x_3 + x_4)^2,
$$
which we called the {\em Descartes form}.
We showed that the set of all Descartes
configurations in the packing formed an orbit of
a discrete subgroup $\sA$ acting on this space,
which was independent of the packing, and
which we termed the {\em Apollonian group.}
This group consisted of  $4 \times 4 $ integer matrices,
and the integrality properties of the Apollonian group
explained the occurrence of packings with
integral curvatures and (curvature)$\times$(circle center)
data. If a single Descartes configuration in
the packing has such integrality properties,
then every Descartes configuration in the
packing inherits such properties, and thus
all the circles in the packing have such properties.

Part I  also introduced two other discrete subgroups
of $Aut(Q_{D})$, which had integrality
properties and a geometric interpretation.
These were the {\em dual Apollonian group} $\sA^{\perp}$
and the group generated by both the Apollonian
group and the dual Apollonian group together,
which we called the {\em super-Apollonian group} $\sA^{S}$. 

In this paper we study to what
extent these results carry over to the $n$-dimensional case.
We show that there
exist $n$-dimensional analogues of the Apollonian
group, dual Apollonian group and super-Apollonian
group. We  study  properties of their orbits
on the parameter space $\sM_{\dd}^n$ of $n$-dimensional
ordered, oriented Descartes configurations  
which we call {\em Apollonian cluster ensembles}.
We determine to what extent properties of the
two-dimensional case extend to the $n$-dimensional
case. One immediate difference is that in dimensions
$n \ge 4$ the orbits of the Apollonian group no
longer correspond to sphere packings; the spheres
in these configurations overlap. However the
groups involved still have rationality properties,
and we consider the question when there exist
ensembles having all curvatures rational and/or
having all data (curvature)$\times$(sphere center)
rational.

In \S2 we describe  six properties that hold in
the two-dimensional case, and then summarize what
the main results of this paper say about generalizations
of these to higher dimensions. Some of them  
generalize completely, others only in 
an infinite set of  specific
dimensions, and some are specific to dimension $2$.

In \S3 we  describe the 
parameter space $\sM_{\dd}^n$ 
of all $n$-dimensional ordered, oriented Descartes configurations
in a coordinate system, ACC-matrices, as those $(n+2) \times (n+2)$
real matrices $\bW$ with $\bW^T \bQ_{D,n} \bW = Q_{W,n}$ where 
$Q_{D,n} = x_1^2 +... + x_{n+2}^2 - \frac{1}{n}(x_1 +... + x_{n+2})^2$ is the
(matrix of the) $n$-dimensional Descartes quadratic form and 
$Q_{W,n} = -8x_1x_2 + 2x_3^2 + ...+ 2x_{n+2}^2$. On the space $\sM_{\dd}^n$
of ACC-matrices the  group
$Aut(Q_{D,n})$ acts on the left, and $Aut(Q_{W,n})$ acts on the right.
Both these  groups are  isomorphic to the 
$(n+2)$-dimensional Lorentz group $O(n+1,1)$, and give two
different ``geometric'' actions. The right action of $Aut(Q_{W,n})$ 
(essentially) corresponds to M\"{o}bius transformations 
acting on the underlying
Euclidean space $\rr^n$ while the left action of $Aut(Q_{D,n})$ is defined
only on the parameter space.

In \S4 we define the $n$-dimensional analogues of the
Apollonian group, the dual Apollonian group and the
super-Apollonian group introduced in part I. As we just
noted above, one immediate
difference with the $2$-dimensional case is that in
dimensions
$n \ge 4$ the orbits of the Apollonian group no
longer correspond to sphere packings; the spheres
in these configurations overlap. Furthermore, even
viewed in the parameter space $\sM_{\DD}^n$ the action of
the $n$-dimensional Apollonian group $\sA_{n}$ is
not discrete. However the action of the dual Apollonian
group $\sA_{n}^{\perp}$ is discrete on $\sM_{\DD}^n$ for all $n \ge 2$.
To restore a discrete group action for the Apollonian
group and the super-Apollonian group in certain higher dimensions, 
it suffices to view it as a diagonal action on
$\sM_{\DD}^n(\RR) \times \prod_{p~|~n-1} \sM_{\DD}^n(\QQ_p)$,
in which $\sM_{\DD}^n(\QQ_p)$ is a $p$-adic parameter space.
This can be done in dimensions for which there are rational
points in the parameter space $\sM_{\DD}^n = \sM_{\DD}^n(\RR)$.

In \S5 we determine presentations for the Apollonian group $\sA_n$
in all dimensions $n \ge 3$. In dimension $3$ there are
extra relations of the form $(\bS_i \bS_j)^3=\bI$
whenever $i \ne j$. The group $\sA_n$ is a hyperbolic
Coxeter group in all dimensions $n \ge 3$.

In \S6 we investigate rationality properties of 
Apollonian cluster ensembles. We prove that strongly
rational configurations exist if and only if the
dimension $n = 2k^2$ or $n = (2k+1)^2$. 

In \S7 we consider to what extent the  duality operator
$\bD$ defined in dimension $2$ has an $n$-dimensional
analogue. In dimension $2$ this operator takes a Descartes
configuration $\sD$ to a dual Descartes configuration $\sD'$
consisting
of the $4$ circles each of which is orthogonal to
all but one of the circles in the packing. We observe
that this operation has a geometric analogue in dimension
$n \ge 3$, which maps the configuration $\sD$ to
the set $\sE$ of  $n+2$ spheres having the property of being  orthogonal
to all but one of the spheres in $\sD$.
These spheres in $\sE$ do not form a Descartes configuration,
but  have an equi-inclination property instead.
In dimensions $n \ge 3$ 
this  operation does not have an algebraic interpretation
in $Aut(Q_{D, n}).$

In \S8 we make some concluding remarks, stating some
unresolved questions.

In the Appendix  we describe the action of the 
$n$-dimensional M\"{o}bius group
M\"{o}b(n) on (ordered, oriented) Descartes configurations, encoded
via an isomorphism, as the action of  $Aut(Q_{W, n})$. 

The general framework of this paper was developed by 
the second author (JCL), who also did much of the writing.
This paper is a revised and extended version of a preprint
originally written in 2000.

\paragraph {Acknowledgments.} 
The authors are grateful for helpful comments from  Andrew Odlyzko
and N. J. A. Sloane. The authors  thank  the reviewer for incisive comments
leading to  reorganization of the paper. 

%
%

\section{Main Results}
\setcounter{equation}{0}
The object of this paper is to study to what extent the 
properties of Apollonian circle packings studied in parts I and II
have $n$-dimensional analogues. 
The two-dimensional case had the following six features: \\

(P1) [Parameter Space Property] The space 
$\sM_{\DD}$ of all ordered, oriented Descartes configurations
can be identified with  the set of real intertwining matrices $\bW$ under 
congruence between two rational
quadratic forms in $4$ variables, the Descartes form $Q_{D}$
and a certain form $Q_{W}$, i.e. $\bW^T \bQ_{D} \bW = \bQ_{W}$.
We called the $4 \times 4$ matrix $\bW$ the augmented curvature-center
coordinates of the associated ordered, oriented Descartes
configuration. The parameter space $\sM_{\DD}$ is a real algebraic
variety and a principal homogeneous
space for the Lorentz group $O(3,1)$.  \\

(P2) [Orbit Property] There exist three groups of matrices
in $ Aut(Q_{D})$ of determinants $\pm 1$, 
 the Apollonian group $\sA$,
the dual Apollonian group $\sA^{\perp}$ and the
super-Apollonian group $\sA^{S}$, which have geometrically
characterizable  actions on Descartes configurations. 
In particular, the
set of Descartes configurations in an Apollonian circle
packing comprise  a single  orbit of a single Descartes
configuration under the action of the Apollonian group
(disregarding  orientation and ordering). \\

(P3) [Integer Matrix Property] The Apollonian group, dual Apollonian group and
super-Apollonian group each consist
of integer matrices. Thus all  three groups are
discrete subgroups of $GL(4, \RR)$.
They each have a discontinuous action on the parameter space $\sM_{\DD}$
of ordered, oriented Descartes configurations. \\

(P4) [Coxeter Group Property] 
 Each of the  Apollonian group, dual Apollonian group, and
super-Apollonian group  have group presentations
as hyperbolic Coxeter groups. \\

(P5) [Integral Packing Property] 
There exist integer Apollonian circle packings,
ones in which all Descartes configurations have
integer curvatures. Furthermore
there exist super-integral Descartes configurations,
ones whose  augmented curvature-center coordinate
matrices $\bW$ are integer matrices. 
There are super-integral Apollonian packings, ones
in which every Descartes configuration has this
property. The set of all
super-integral ordered, oriented
Descartes configurations forms $672$ orbits under
the action of the super-Apollonian group. \\

(P6) [Duality Operator]
There exists a duality operator $\bD$ in $Aut(Q_{\sD})$ given
by a fixed matrix  $GL(4, \RR)$ with half-integer entries.  
It acts by conjugacy on the super-Apollonian group and  gives an
outer automorphism of order $2$ of this group.
 This outer automorphism   
conjugates the Apollonian group to the dual
Apollonian group. The duality operator has a
geometric interpretation in terms of its action
on individual Descartes configurations. \\

The results of this paper generalize each of these
features, which we consider in order.
The parameter space property (P1)
generalizes to all dimensions, as was
shown by three of the authors in \cite{LMW02}. Results in
this direction were established
earlier by  Wilker \cite{Wi81}, who used
the term {\em cluster} for an ordered, but not
oriented, Descartes configuration.
There is a notion  of augmented curvature-center
coordinates (ACC-coordinates) for 
an ordered, oriented
$n$-dimensional  Descartes configuration, and a parameter
space $\sM_{\DD}^n$ of all 
such  configurations, specified by a matrix condition that 
intertwines two quadratic forms in $n+2$ variables
under conjugacy, as was shown in \cite{LMW02}.
Furthermore the space $\sM_{\DD}^n$ has the
structure of a principal homogeneous
space for $O(n+1, 1)$. We present these results in \S3.

The orbit property (P2) generalizes to all dimensions. There exists 
an $n$-dimensional version of the Apollonian group,
the dual Apollonian group and the super-Apollonian
group, with the same geometric interpretation of
their actions on Descartes configurations. 
These are described in \S4. We define an 
{\em Apollonian cluster ensemble} to be an orbit of
a single (ordered, oriented) Descartes configuration under the action
of the Apollonian group. However there is a sphere
packing interpretation of the geometric packing
corresponding to this orbit in dimension
$3$ only, and not in dimensions $n \ge 4$.

The integer matrix  property  (P3) partially generalizes to all dimensions,
as we discuss in \S4. 
The $n$-dimensional Apollonian group,
dual Apollonian group and super-Apollonian
group each  consist of  integral matrices in dimension
$3$, and they act discretely. 
The dual Apollonian
group consists of integer matrices in all dimensions,
and so is a discrete subgroup of $GL(n+2, \RR)$ and
acts discontinuously on the parameter space $\sM_{\DD}^n$.
However for $n \ge 4$ the Apollonian group and super-Apollonian
group consist of rational matrices, with denominators whose
prime factors all divide $n-1$. The Apollonian group $\sA_n$
is not a discrete subgroup of $GL(n+1, \RR)$ for $n \ge 4$
and does not  have a discontinuous action on $\sM_{\DD}^n$
(Theorem~\ref{le4.1.1}).
However we can restore discreteness by adding an action on
some $p$-adic groups, for those primes $p$ dividing $n-1$.
That is, the super-Apollonian
group $\sA_n^{S}$ can be embedded discretely by a diagonal
action inside $GL(n+2, \RR) \times \prod_{p~|~n-1} GL(n+2, \QQ_p)$
for all $n$. For certain dimensions $n$, 
those characterized in \S6 as $n=2k^2$ or
$n=(2k+1)^2$ for $k \ge 1$, 
one can get a discontinuous diagonal action of the 
super-Apollonian group on a parameter space
$\sM_{\DD}^n(\RR) \times \prod_{p~|~n-1}  \sM_{\DD}^n(\QQ_p).$
For this it is necessary that the p-adic parameter spaces
$\sM_{\DD}^n(\QQ_p)$ have enough points. Finally on
the level of ``packings'' in the weak sense of spheres
not crossing each other, this holds in dimensions $2$
and $3$ for all these groups, and holds for the dual
Apollonian group $\sA_n^{\perp}$ for all $n \ge 4$ (Theorem \ref{le4.1.3}),
but not for the other two groups.

The hyperbolic Coxeter group property (P4) 
partially (and perhaps completely) generalizes to all dimensions.
We show in \S5 that for all $n \ge 3$
the $n$-dimensional Apollonian
group is a hyperbolic Coxeter group
(Theorems ~\ref{th51} and \ref{th52}). 
For $n=3$ the group has an extra
relation, which explains the existence of the structures
``The Bowl of Integers'' and ``The Hexlet'' studied
by Soddy \cite{Sod36b}, \cite{Sod37a} \cite{Sod37b} and Gosset \cite{Go37b}. 
It may be true  that in all dimensions the dual
Apollonian group and super-Apollonian group
are also   hyperbolic Coxeter groups, but we
leave these as open questions.

The integral packing properties in (P5) 
generalize to all dimensions in weakened forms.
In \S6 we observe that in all dimensions $n \ge 3$ there exist
Apollonian cluster ensembles in which every sphere
has a rational curvature, with denominators
divisible only by a certain finite set of primes, those
dividing $n-1$ (Theorem~\ref{th61}). We
then consider the question in which dimensions $n$ does
there exist a super-rational Descartes
configuration, i.e. one whose augmented curvature-center
coordinate matrix $\bW_{\sD}$ is 
a rational matrix in $\sM_{\DD}^n$.
We prove this can be done if and only if
the dimension $n = 2k^2$ or $n = (2k+1)^2$ for 
some $k \ge 1$ (Theorem \ref{th63}).

The duality operation in (P6) giving a conjugacy between
the Apollonian group and the dual Apollonian
group, appears to exist algebraically  only in dimension 2.
In \S7 we show that viewed as a  geometric action, 
there does exist a natural ``duality operator''
acting on Descartes configurations in all dimensions $n \ge 2$.
However in dimensions $n \ge 3$, the image of this operator is
not a Descartes configuration, but 
instead is a collection of $n+2$ spheres
each of which intersects the others in a particular fixed angle
$\theta_n$, which depends
on the dimension $n$ (Theorem~\ref{th72}).

%
%
%

\section{Descartes Configurations and Group Actions in $\RR^n$}
\setcounter{equation}{0}

%
%

\subsection{Descartes Configurations and ACC-Coordinates}  

We start with the
 generalization of the Descartes circle theorem to the
$n$-dimensional case. An {\em Descartes configuration} in $\RR^n$
consists of $(n+2)$ pairwise tangent $(n-1)$-spheres 
$(S_1, S_2, \cdots, S_{n+2})$ in $\RR^n$,
with all points of tangency distinct. 
This result, which was  
termed the {\em Soddy-Gossett theorem} in \cite{LMW02}, 
after Soddy \cite{Sod36a}     and Gossett \cite{Go37} 
(although it was discovered earlier in the $3$-dimensional case), 
states that
if the spheres have disjoint interiors then 
\beql{003a}
\sum_{j=1}^{n+2} \frac{1}{r_i^2} = 
\frac{1}{n} \left(\sum_{i=1}^{n+2} \frac{1}{r_i}\right)^2.
\eeq
The Descartes circle theorem is the case $n=2$.

The Soddy-Gossett theorem
holds for  all Descartes configurations,
including configurations where one sphere encloses the others,  
provided that we assign appropriate signs to the curvatures,
so that the configuration has a (total) orientation, as defined below.
An {\em oriented sphere} is a sphere
together with an assigned direction of unit normal vector, which can
point inward or outward.  If it has radius $r$ then its {\em oriented
radius} is $r$ for an inward pointing normal and $-r$ for an outward
pointing normal.  Its {\em oriented curvature} (or ``signed
curvature'') is $1/r$ for an inward pointing normal and $-1/r$ for an
outward pointing normal.  By convention, the {\em interior} of an
oriented sphere is its interior for an inward pointing normal and its
exterior for an outward pointing normal.  An {\em oriented Descartes
configuration} is a Descartes configuration in which the orientations
of the spheres are compatible in the following sense:  either (i) the
interiors of all $n+2$ oriented spheres  are disjoint, or (ii) the
interiors are disjoint when all orientations are reversed.  Each
Descartes configuration has exactly two compatible orientations in
this sense, one obtained from the other by reversing all orientations.
The {\em positive (total)  orientation} of a Descartes configuration is
the one in which the sum of the signed curvatures is positive, while
the {\em negative (total)  orientation} is the one in which the sum of
the curvatures is negative. One can check that the sum of the 
curvatures cannot be zero. 

Now let  $b_i= \frac{1}{r_i}$ denote the (signed) curvature of
the $i$-th sphere of an ordered, (totally) oriented Descartes 
configuration , and  let $\bb= (b_1, ..., b_{n+2})$. 
The geometry of such a 
 Descartes configurations is encoded in the signed curvature
vector $\bb.$ In
the positively oriented case, where  $\sum_{j = 1}^{n+2} b_j > 0$, 
one of the
following holds:  (i) all of $b_1, ~b_2,\ldots,~b_{n+2}$ are positive;
(ii) $n+1$ are positive and one is negative; (iii) $n+1$ are positive
and one is zero; or (iv) $n$ are positive and equal and the other two
are zero.  These four cases correspond respectively to the following
configurations of mutually tangent spheres:  (i) $n + 1$ spheres, with
another in the curvilinear simplex that they enclose; (ii) $n+1$
spheres inscribed inside another larger sphere; (iii) $n+1$ spheres
with one hyperplane (the $(n + 2)$-nd ``sphere''), tangent to each of
them; (iv) $n$ equal spheres with two common parallel tangent planes.

We can reformulate the Soddy-Gossett theorem in matrix terms.
Let $b_i= \frac{1}{r_i}$ denote the (signed) curvature of
the $i$-th sphere of an ordered, oriented Descartes 
configuration , and  let
 $\bb= (b_1, ..., b_{n+2})$, then \eqn{003a} becomes
\beql{304} 
\bb^T \bQ_{D,n} \bb = 0,
\eeq
in which $\bQ_{D,n}$ is the symmetric matrix associated
to the $n$-dimensional
Descartes quadratic form, defined below.

\begin{defi}~\label{de31} 
{\em
The {\em $n$-dimensional Descartes quadratic form} $Q_{D,n}$
is the quadratic form in $n+2$ variables whose
associated symmetric matrix is
\beql{403a}
\bQ_{D, n} := \bI_{n+2} - \frac{1}{n} {\bf 1}_{n+2} {\bf 1}_{n+2}^T,
\eeq
in which ${\bf 1} = (1, 1, ... , 1)^T$ is a column vector of
length $n+2$.
}
\end{defi}

Here the original Descartes quadratic form is 
\beql{003}
\bQ_{D,2}= I - \frac{1}{2}{\bf 1}{\bf 1}^T =          
         \frac{1}{2}\left[ \begin{array}{rrrr}
               1 & -1 & -1 & -1 \\
               -1 & 1 & -1 & -1 \\
               -1 & -1 & 1 & -1 \\
               -1 & -1 & -1 & 1 \end{array} \right].   
\eeq

In \cite{LMW02} three of the authors of this paper 
showed  that there exists
a parametrization of the set $\sM_{\DD}^n$ of all ordered, 
oriented Descartes configurations
in  $\RR^n$,
using  a coordinate system involving the
curvatures and centers of the spheres, which appears
as Theorem~\ref{th31} below.

\begin{defi}~\label{de32}
{\em
Given an oriented sphere $S$ in $\rr^n$, its
{\em curvature-center coordinates} consist of the $(n+1)$-vector
\beql{300d}
        \bm(S) = (b, b x_1, \ldots, b x_n),
\eeq
in which $b$ is the signed curvature of $S$ (assumed nonzero) and
$\bx(S) = \bx = (x_1, x_2, \ldots, x_n)$ is its center.  For the
degenerate case of an oriented hyperplane $H$, its 
{\em curvature-center coordinates} $\bm(H)$ are defined to be
\beql{300e}
        \bm(S) = (0,\bh),
\eeq
where $\bh :=  (h_1, h_2, \ldots, h_n)$ is the unit normal vector
that gives the orientation of the hyperplane.
}
\end{defi}

To see the origin of this definition in the degenerate case, let the
point of $H$ closest to the origin be $\bz = \lambda \bh$ for some real
value $\lambda$.  For $t > |\lambda|,$ let $S_t$ be the oriented sphere 
of radius
$t$ centered at $(t + \lambda)\bh,$ which has center in direction $\bh$ from
the origin and contains $\bz$.  As $t \to \infty$ the oriented spheres
$S_t$ clearly converge geometrically to the oriented hyperplane $H$, and
$\bm(S_t) = (\frac{1}{t}, (1 + \frac{\lambda}{t})\bh) \to \bm(H) = (0, \bh).$

Curvature-center coordinates are not quite a global coordinate system,
because they do not always uniquely specify an oriented sphere.  Given
$\bm \in \rr^{n+1}$, if its first coordinate $b$ is nonzero then there
exists a unique sphere having $\bm = \bm(S)$.  But if $b = 0$, the
hyperplane case, there is a hyperplane if and only if $\sum h_i^2 =
1,$ and in that case there is a pencil of hyperplanes that have the
given value $\bm$, which differ from each other by a translation.

We obtain a global coordinate system for spheres by
adding an additional coordinate. This coordinate incorporates
information about the sphere $\bar{S}$ obtained from $S$ by
inversion in the unit sphere. 
In $n$-dimensional Euclidean space, the operation of {\em inversion in
the unit sphere} replaces the point $\bx$ by $\bx/{{|\bx|}^2}$, where
$|\bx|^2 = \sum_{j = 1}^n x_j^2$.  Consider a general oriented sphere
$S$ with center $\bx$ and oriented radius $r$.  Then inversion in the
unit sphere takes $S$ to the sphere $\bar{S}$ with center $\bar{\bx} =
{\bx}/(|\bx|^2 - r^2 )$ and signed radius $\bar{r} = r/(|\bx|^2 -
r^2 )$.  If ${|\bx|}^2 > r^2$, then $\bar{S}$ has the same sign
as $S$.  In all cases,
\beql{303f}
        \frac {\bx}{r} = \frac{\bar{\bx}}{\bar{r}}, 
\eeq
and
\beql{306aa}
        \bar{b} = \frac{|\bx|^2}{r}- r.
\eeq

\begin{defi}~\label{de33}
{\em
Given an oriented sphere $S$ in $\rr^n$, its
{\em augmented curvature-center coordinates} (or {\em ACC-coordinates}) 
of $S$ are 
given by the $(n+2)$-vector
\beql{300g}
        \bw(S) := ( \bar{b}, b, b x_1, \ldots , b x_n) = (\bar{b}, \bm),
\eeq
in which $\bar{b}= b(\bar{S})$ is the curvature of the sphere or
hyperplane $\bar{S}$ obtained by inversion of $S$ in the unit sphere,
and the entries of $\bm$ are its curvature-center coordinates.  For
hyperplanes we define
\beql{300f}
        \bw(H) := (\bar{b}, 0, h_1, \ldots , h_n) = (\bar{b}, \bm),
\eeq
where $\bar{b}$ is the oriented curvature of the sphere or hyperplane
$\bar{H}$ obtained by inversion of $H$ in the unit sphere.
}
\end{defi}

Augmented curvature-center coordinates provide a global coordinate
system for oriented spheres:  
no two distinct oriented spheres have the same coordinates.
The only case to resolve is when $S$ is a hyperplane, i.e., $b=0$.
The relation \eqn{303f} shows that $(\bar{b}, b x_1, \ldots  , b x_n)$
are the curvature-center coordinates of $\bar{S}$, and if $\bar{b}
\neq 0$, this uniquely determines $\bar{S}$; inversion in the unit
circle then determines $S$.  In the remaining case, $b = \bar{b} =0$
and $S = \bar{S}$ is the unique hyperplane passing through the origin
whose unit normal is given by the remaining coordinates.

Given a collection $(S_1, S_2, \ldots  , S_{n + 2})$ of $n + 2$ oriented
spheres (possibly hyperplanes) in $\rr^n$, the {\em augmented matrix}
$\bW$ associated with it is the $(n+2) \times (n + 2)$ matrix whose
$j$-th row has entries given by the augmented curvature-center
coordinates $\bw(S_j)$ of the $j$-th sphere.


To state the next result, we introduce another quadratic form.

\begin{defi}~\label{de34}
{\em 
The {\em ($n$-dimensional) Wilker quadratic form} $Q_{W, n}$
 is the $(n+2)$-variable
quadratic form given by the symmetric $(n+2) \times (n+2) $ matrix
\beql{307a}
\bQ_{W, n} : =
        \left[ \begin{array}{rrl}
           0 & -4 & 0   \\
          -4 &  0 & 0   \\
           0 &  0 & 2\bI_n
        \end{array} \right].
\eeq
}
\end{defi}

This name is made in honor of  J. B. Wilker ~\cite{Wi81}, who
introduced in spherical geometry a coordinate system
analogous to augmented curvature-center coordinates, see
\cite[\S2 p. 388-390 and \S9]{Wi81}. However he did not
formulate any result explicitly exhibiting a quadratic form
like $\bQ_{W, n}$; see the remark on \cite[p. 349 ]{LMW02}.

%
%
%

\begin{theorem}~\label{th31}
{\em [Augmented Euclidean Descartes Theorem]}
The augmented matrix $\bW = \bW_{\sD}$ of an oriented Descartes configuration
$\sD$ of $n + 2$ spheres $\{S_i:  1 \leq i \leq n+2\}$ in $\rr^n$ satisfies
\beql{307}
        \bW^T \bQ_{D,n}  \bW = 
        \left[ \begin{array}{rrl}
           0 & -4 & 0   \\
          -4 &  0 & 0   \\
           0 &  0 & 2\bI_n
        \end{array} \right].
\eeq
Conversely, any real solution $\bW$ to \eqn{307} is the augmented
matrix $\bW_{\sD}$ of a unique ordered, oriented Descartes configuration
$\sD$.  
\end{theorem}

\paragraph{Proof.} 
This is proved as Theorem 3.3  in \cite{LMW02}. $~~~\bsq$ \\

Theorem~\ref{th31} 
states that the  augmented curvature-center coordinates
of an ordered, oriented Descartes configuration give an intertwining map
between the Descartes form and the Wilker form. 
We note that the Soddy-Gossett theorem is a 
special case of Theorem~\ref{th31},
encoded as the  $(2,2)$ entry of the matrix $\bW^T \bQ_{D,n}  \bW$.

Both the Descartes quadratic form and the Wilker quadratic form
are equivalent over the real numbers to the Lorentzian
quadratic form 
\beql{313aa}
Q_{\sL, n}(x) := - x_0^2 + x_1^2 + \cdots + x_{n+1}^2,
\eeq
see \S3.2. 
This quadratic form has a large group of real automorphisms
under congruence 
$$
Aut(Q_{\sL, n}) = \{ \bU  \in GL(n, \RR) : 
\bU^T \bQ_{\sL, n} \bU = \bQ_{\sL, n} \},
$$
which is the Lorentz group $O(n+1,1)$.
In consequence both the Descartes quadratic form $Q_{D,n}$
and Wilker quadratic form $Q_{W,n}$ have automorphism groups under
(real) congruence which are conjugates of the Lorentz
group.

The Descartes form is not only equivalent
to the Lorentz form  over the real numbers,
but sometimes  over the rational numbers. 
In dimension $2$, the Descartes and Lorentz forms 
 are rationally equivalent, where one has 
\beql{N313}
\bQ_{\sL, 2}(x) = \bJ_0^T \bQ_{D, 2} \bJ_0,
\eeq
with 
$$
\bJ_0 = \frac{1}{2} \left[
\begin{array}{crrr}
1 & 1 & 1 & 1 \\
1 & 1 & -1 & -1 \\
1 & -1 & 1 & -1 \\
1 & -1 & -1 & 1
\end{array}
\right]
$$
In dimension $2$ the Wilker form $Q_{W,2}$ is also rationally
equivalent to both the Descartes form and Lorentz form,
as shown in part I, \S3.1.  

In higher dimensions these three forms are not
always rationally equivalent.
A necessary condition for rational equivalence of two
quadratic forms is that their determinants differ
by a rational square. We have $\det(\bQ_{\sL}) = -1$,
$\det(\bQ_{D}) = -\frac{2}{n}$ and 
$\det(\bQ_{W}) = -2^{n+4}$. It follows that a 
necessary condition for rational equivalence of
the Descartes and Lorentz forms in dimension $n$ is
that $n= 2k^2$;  for the Wilker and Lorentz forms in dimension $n$
that $n = 2k$; and for the Descartes and Wilker forms
in dimension $n$ that  $n= k^2$ for odd $k$ or $n=2k^2$ for even $k$.
All three of these necessary conditions hold if and only if
$n=2k^2$.

In \S6 we show that the last condition is sufficient
for equivalence of the Descartes and Wilker forms.

%
%

\subsection{M\"{o}bius and Lorentz group actions}

The augmented Euclidean Descartes theorem  
yields two group actions on the space of Descartes configurations.
The group $Aut( Q_{D, n})$ acts on the left
and the group $Aut(Q_{W, n})$ acts on the right, as
$$
\bW_{\sD}  \mapsto \bU  \bW_{\sD}  \bV^{-1}, ~~\mbox{with}~~ 
\bU \in Aut( Q_{D, n}),~ 
\bV \in Aut(Q_{W, n}).
$$
The two group actions obviously commute with each other.
The following result generalizes the two-dimensional case,
where the Lorentzian form $Q_{\sL, 2}$ is given in \eqn{N313}.

%
%
%

\begin{theorem}~\label{nth32}
(1) The groups $Aut(Q_{D,n})$ and $Aut(Q_{W,n})$ are
conjugate over the real numbers to $Aut(Q_{\sL,n}) \equiv  O(n+1,1)$.

(2) The group $Aut(Q_{D,n})$ acts transitively on the 
left on the parameter space $\sM_{\DD}^n$
of all ordered, oriented  Descartes configurations. 
Given two such 
 Descartes configurations $\sD$ and $\sD^{'}$
there exists a unique $\bU \in Aut(Q_{D,n})$ such that
$\bU \bW_{\sD} = \bW_{\sD'}$.

(3) The group $Aut(Q_{W})$ acts transitively on the 
right on the space 
of all ordered, oriented  Descartes configurations $\sM_{\DD}^n$. 
Given two such Descartes configurations $\sD$ and $\sD^{'}$
there exists a unique $\bV \in Aut(Q_{W,n})$ such that
$\bW_{\sD} \bV^{-1} = \bW_{\sD'}$.
\end{theorem}

\paragraph{Remark.} This result allows one to define both
a left and right  $O(n+1,1)$ action 
on the space $\sM_{\DD}^n$, which  depends on the
choice of conjugacy made in (1). Then (2) and (3) show 
both these actions are transitive and
have trivial stabilizer.  This gives $\sM_{\DD}^n$
the structure of a principal homogeoneous space 
(or torsor) for $O(n+1,1)$, for each of these
actions. 

\paragraph{Proof.}
Part (1) follows for $Q_{W,n}$ on taking
$\bQ_{\sL, n} = \bZ^T \bQ_{W,n} \bZ$
with 
$$
\bZ = \frac{1}{\sqrt{2}}\left[ \begin{array}{rrr}
\frac{1}{2} & -\frac{1}{2} & 0 \\
\frac{1}{2}& \frac{1}{2} & 0 \\
0 & 0 & {\bf I}_{n}
\end{array} 
\right] .
$$
It then follows for $Q_{D,n}$ because it is conjugate
to $Q_{W,n}$ by Theorem~\ref{th31}, using $\bW_{\sD}$
for any fixed Descartes configuration.

Parts (2) and (3) follow immediately from (1),
because $Aut(Q_{\sL,n}) = O(n+1, 1)$ acts transitively.
~~~$\bsq$ \\

Since the spheres in a Descartes configuration appear
as the rows in the augmented matrix $\bW_{\sD}$ of an oriented
Descartes configuration, 
the action on the right by elements of $Aut(Q_{W, n})$
maps spheres to spheres. This action 
can essentially be identified
with the M\"{o}bius group of linear fractional
transformations acting on the one-point compactification
$\hat{\RR}^n$ of $\RR^n$. More precisely, it
corresponds to a direct product of the M\"{o}bius group
with $\{ \bI, -\bI\}$,  because the M\"{o}bius group preserves
total orientation of Descartes
configurations. A precise description of the isomorphism
is given in the Appendix. (The case $n=2$ was treated
in Appendix A of part I.)

The action on the left, by $Aut(Q_{D, n})$,
mixes together the 
different spheres in the original   Descartes configuration,
and does not make sense as an action on individual spheres. 
This group action is intrinsically associated to
the $\frac{n(n+1)}{2}$-dimensional (real) 
parameter space $\sM_{\dd}^n$ of oriented
Descartes configurations.

%
%
%

\section{Apollonian Groups and Apollonian Cluster Ensembles}
\setcounter{equation}{0}
In parts I and II we studied 
Apollonian circle packings 
in terms of the Descartes configurations they
contain. We showed  they consisted of a single
orbit of a discrete subgroup of the automorphism
group $Aut(Q_{D, 2})$ of the Descartes quadratic
form $Q_{D, 2}$. We also introduced there
another discrete subgroup of $Aut(Q_{D, 2})$,
the dual Apollonian group, and in addition,
the super-Apollonian group, which is the subgroup
of $Aut(Q_{D, 2})$ generated by the Apollonian
group and the dual Apollonian group together.

We show there are 
analogues of these three groups in all dimensions
$n \ge 3$, which are subgroups of $Aut(Q_{D, n})$
consisting of rational matrices. 
We call these the $n$-dimensional Apollonian group,
dual Apollonian group, and super-Apollonian group.
We then define an {\em Apollonian cluster ensemble} to be
an orbit of the Apollonian group. This provides a
generalization of Apollonian packing to all dimensions,
though it turns out not to correspond to a sphere-packing
in dimensions $n \ge 4$.

%
%
%
\subsection{$n$-Dimensional Apollonian Group}

The $n$-dimensional Apollonian group 
$\sA_n = \langle \bS_1, \bS_2, ..., \bS_{n+2} \rangle$ 
consists of $(n+2) \times (n+2)$ matrices with 
\beql{N401a}
\bS_j := \bI + \frac{2}{n-1} \be_j \bo^T - \frac{2n}{n-1} \be_j \be_j^T,
\eeq
where $\be^j$ is the $j$-th unit coordinate (column) vector and 
$\bo= \be_1 +\cdots +\be_{n+2} =(1,1, ..., 1)^T$.
That is, 
$\bS_j$ is  the identity matrix in
all rows but the $j$-th row, and there has $-1$ on
the diagonal, and all off diagonal elements equal
to $\frac{2}{n-1}$. It is straightforward to check that 
$\sA_n \subset \mbox{Aut}(Q_{D,n})$. The relations
$\bS_j^2 = \bI$ are evident. 

The algebraic action of the operator $\bS_j$ on 
the augmented curvature center coordinates 
$\bW_{\sD}$ of
a Descartes configuration $\sD$
 is to take it to $\bS_j\bW_{\sD}$,
which is  $\bW_{\sD'}$ for some $\sD'$. 
The geometric interpretation of this action 
is to  fix all (oriented) spheres in the Descartes configuration $\sD$
except the $j$-th sphere, and to replace that sphere with the
unique other sphere that is tangent to the remaining
$n+1$ spheres, assigned an appropriate orientation; this is $\sD'$.
For a fixed Descartes
configuration $\sD$  this  operation can also be realized by a
M\"{o}bius transformation that is an inversion with
respect to the sphere that passes through the $\frac{n(n+1)}{2}$
tangency points of the remaining $n+1$ spheres in the
Descartes configuration. (The existence of such a sphere
is demonstrated in Theorem~\ref{th71}.)
  
It is apparent that $\sA_n$ is a group of integer matrices
for $n= 2$ and $3$, while for $n \ge 4$ it  is not always
integral, consisting of rational
matrices whose denominators contain only powers of
primes that divide $n-1$.  
%
%
%

\begin{theorem}~\label{le4.1.1}
The Apollonian group $\sA_n$ is a
discrete subgroup of $GL(n+2, \RR)$ for $n=2$ and $n=3$.
It is not a discrete subgroup
of $GL(n+2, \RR)$ for all $n \ge 4.$
\end{theorem} 

\paragraph{Proof.}
The group $\sA_n$ is discrete in dimensions $2$ and $3$
since it is then a subgroup of \\
$GL(n+1, \ZZ).$ 

We consider the element $\bS_1 \bS_2$, and show for $n \ge 4$
that the  set $\{ (\bS_1 \bS_2)^k : k \ge 1$ does not contain the
identity matrix, but the closure of this set does contain
the identity matrix. The $(n+1) \times (n+1)$
matrix $\bS_1 \bS_2$ is a product
of two reflections, so it has
determinant $1$ and has $n$ of its eigenvalues equal to $1$.
Its first two rows are $(a_n^2-1,  -a_n, a_n^2+a_n, a_n^2+a_n, ...)$
and $(a_n, -1, a_n, a_n, ...)$ in which $a_n = \frac{2}{n-1}$,
and all other rows are those of the identity matrix.
Its two non-unit  
eigenvalues are $e^{\pm i \theta_n}$ with
$\cos \frac{\theta_n}{2} = \frac{a_n}{2}$, so 
$\theta_n = 2 \cos^{-1}(\frac{1}{n-1})$. Since these
two eigenvalues are distinct, it  follows that
$\bS_1 \bS_2$ is diagonalizable. Now 
$\Sigma_n := \{(\bS_1 \bS_2)^k:~ k \ge 1\}$ 
contains diagonalizable elements with all
eigenvalues arbitrarily close to $1$, so it contains the
identity matrix as a limit point. When $n \ge 4$ it is
well-known that
$\theta_n$ is an irrational multiple of $\pi$, 
whence  $\Sigma_n$ does not contain the identity matrix,
which shows that the Apollonian group $\sA_n$ is not
a discrete subgroup of $GL(n+1, \RR)$.
~~~$\bsq$ \\

Because the Apollonian group is arithmetic 
one can view it as a discrete group if one considers its
actions on certain 
$p$-adic groups $GL(n+2, \QQ_p)$ for $p | n-1$. 
Let $n^{\ast} = n-1$ if $n$ is even, and $\frac{n-1}{2}$
if $n$ is odd. 
Then one can establish that  the Apollonian group $\sA_n$
is a discrete group when embedded diagonally
in $GL(n+2, \RR) \times \prod_{p~ |~n^{\ast}}GL(n+2, \QQ_p)$.
Using this structure, one may
obtain a  discrete
(i.e. discontinuous) diagonal action of the Apollonian group 
on the space
$\sM_{\DD}^n(\RR) \times  \prod_{p~ |~n^{\ast}} \sM_{\DD}^n(\QQ_p)$,
whenever all of the parameter
 spaces $\sM_{\DD}^n(\QQ_p)$ are nonempty.
Here the  spaces $\sM_{\DD}^n(\QQ_p)$ are defined as the $p$-adic
solutions to the conditions in \eqn{307}.
The criterion  in \S6.2, shows that $\sM_{\DD}^n(\RR)$
contains rational solutions, hence $p$-adic solutions
for all $p$, when  $n=2k^2$ 
or $(2k+1)^2$, for some $k \ge 1$.
%
%
%
\subsection{$n$-Dimensional Dual Apollonian Group}

The {\em $n$-dimensional dual Apollonian group}
$\sA_n^{\perp}= 
\langle \bS_1^{\perp}, \bS_2^{\perp},..., \bS_{n+2}^{\perp} \rangle$
consists of $(n+2) \times (n+2)$ matrices with 
\beql{N402a}
\bS_j^{\perp} := \bI + 2 \bo \be_j^T - 4 \be_j \be_j^T.
\eeq
That is, 
$\bS_j^{\perp}$ is the identity matrix in all columns except
the $j$-th column, where it has $-1$ as a diagonal entry
and $2$ for all off-diagonal entries. The group $\sA_n^{\perp}$
is a group of integer matrices in all dimensions $n$.
It is straightforward to check that 
$\sA_n^{\perp} \subset \mbox{Aut}(Q_{D,n}).$ It is
evident that $(\bS_j^{\perp})^2 = I$ holds for all $j$.
The geometric interpretation of the operation $\bS_j^{\perp}$ 
on a Descartes configuration $\sD$ is that it
encodes inversion with respect to
 the $i$-th sphere of that  configuration.

%
%
%

\begin{theorem}~\label{le4.1.2}
The dual Apollonian group $\sA_n^{\perp}$ is a
discrete subgroup of $GL(n+2, \RR)$ for all $n \ge 2$.
\end{theorem} 
\paragraph{Proof.}
This holds since  $\sA_n^{\perp}$ is a subgroup of 
the discrete group $GL(n+2, \ZZ)$.
$~~~\bsq$ \\

Orbits of this group, acting on Descartes configurations,
retain a ``packing'' property at the level of
individual spheres.
%
%
%

\begin{theorem}~\label{le4.1.3}
For all $n \ge 2$,  an orbit of the dual Apollonian group 
$\sA_n^{\perp}$ acting on
a single Descartes configuration
gives a  ``packing'' of spheres in the
weak sense that no  two  spheres in distinct Descartes
configurations of the orbit cross each other, i.e. any two such
spheres either coincide, or are disjoint, or are tangent.
\end{theorem} 

\paragraph{Proof.}
This can be seen geometrically by constructing the dual 
packing starting from a single Descartes configuration $\sD$.
The spheres in the resulting ``packing'' are all nested
inside spheres of $\sD$, and we call the level of a sphere
its depth of nesting inside some  sphere of $\sD$,
the spheres in $\sD$ being assigned level $0$.
One proceeds in stages, where at stage $k-1$ one has the
set of Descartes
configurations $\{ \bS_{i_{k-1}}^{\perp} \bS_{i_{k-2}}^{\perp} \cdots 
\bS_{i_1}^{\perp} \bW_{\sD}\}$,
in which each  $i_j \ne i_{j-1}$.
We assert as  an induction hypothesis, that 
 each stage $k-1$ Descartes configuration has all
but one of its spheres at level $k-1$, with one
sphere at level $k-2$. This holds for the base case $k=2$
by inspection. 
At stage $k$ each  multiplication by a generator takes a particular Descartes
configuration at stage $k-1$ 
and maps $(n+1)$ of its  spheres inside
one level $k-1$ sphere of the configuration. 
(Here the condition $i_k \ne i_{k-1}$ is used.)
The new Descartes configuration
then consists of $(n+1)$ spheres nested to depth $k$ inside
$\sD$, contained in  one outer sphere nested to depth $k-1$. 
Furthermore each  level  $k$ Descartes configuration is nested inside
a unique level $k-1$-sphere, of which there are
$(n+2) \cdot (n+1)^{k-2}$ choices.  It follows that
the new level $k$ spheres cannot cross any spheres at any 
levels up to $k-1$,  
and they also cannot cross any other level $k$ spheres 
because they are either inside different level $k-1$ spheres,
or if they are in the same level
$k-1$ sphere, then they form part of  a single Descartes configuration.
Thus the non-crossing property holds at depth $k$,
and the result follows by induction on $k$.$~~~\bsq$

%
%
%
\subsection{$n$-Dimensional  Super-Apollonian Group}

The {\em $n$-dimensional super-Apollonian group} 
$\sA_n^{S}$ is the group generated by $\sA_n$ and $\sA_n^{\perp}$,
so that
$$
\sA_n^{\perp} = \langle \bS_1, \bS_2, \cdots \bS_{n+2}, 
\bS_1^{\perp}, \bS_2^{\perp}, \cdots, \bS_{n+2}^{\perp}\rangle.
$$
This group consists of integer matrices when $n=2$ or $3$, and
of rational matrices otherwise. In particular $\sA_n^{S}$ is
a discrete subgroup of $GL(n+2, \RR)$ for $n=2$ or $n=3$, and
is not a discrete subgroup for $n \ge 4$. 
The property that some spheres in different Descartes
configurations intersect non-tangentially when $n \ge 4$
is inherited from the action of the Apollonian group.

We do not know if the super-Apollonian group is a
Coxeter group. We do know that its generators satisfy the 
Coxeter relations given in the following lemma.
 It seems plausible that for $n \ge 4$ these
are a generating set of relations; if so, then $\sA_n^{S}$ would be
a hyperbolic Coxeter group.
 
\begin{theorem}~\label{le4.2.1}
The super-Apollonian group 
$\sA_n^{S}=\langle \bS_1, \bS_2, \cdots \bS_{n+2}, 
\bS_1^{\perp}, \bS_2^{\perp}, ..., \bS_{n+2}^{\perp}\rangle $
has generators satisfying the relations.
\beql{421}
\bS_j^2 = I ~~\mbox{and}~~~ (\bS_j^{\perp})^2 =I, ~~~ 1 \le j \le n+2.
\eeq
and 
\beql{422}
\bS_j \bS_k^{\perp} = \bS_k^{\perp} \bS_j ~~~\mbox{when}~~ j \ne k.
\eeq
\end{theorem}
\pf
This is a direct calculation from the
definition.  To appreciate \eqn{422} one  can consider more generally
$(n+2) \times (n+2)$ matrices $\bS_j(\lambda)$ 
whose entries are the identity matrix
except in the $j$-th row, where they are $-1$ on the diagonal
and $\lambda$ off the diagonal. Then $\bS_j(\lambda)^2 = I$
holds for $1 \le j \le n$, but the extra relations \eqn{422}
hold only when $\lambda = \frac{2}{n-1}$. (Compare the 
 $(j, k)$-th entry of the products
$\bS_j\bS_k^{\perp}$ and $ \bS_k^{\perp} S_j$). Note however that
$ \bS_j(\lambda) \in \mbox{Aut}(Q_{D,n})$ 
if and only if  $\lambda = \frac{2}{n-1}$.
~~~$\bsq$ \\

As mentioned above for the Apollonian group, there is
a discrete action of the super-Apollonian group
in certain dimensions, 
given by its  diagonal action of $\sA_{n}^{S}$ on 
the product parameter  space
$\sM_{\DD}^n(\RR) \times  \prod_{p~ |~n^{\ast}} \sM_{\DD}^n(\QQ_p)$,
provided that all  the  spaces $\sM_{\DD}^n(\QQ_p)$ are  nonempty.
As shown in \S6.2, this will be the case when $n = 2k^2$
or $(2k + 1)^2$, for some $k \ge 1$.

%
%
%
\subsection{Apollonian Cluster Ensembles}

In the two-dimensional case Apollonian packings could be
described as (a) collections of circles, or (b) the orbit of
four circles under a certain discrete group of M\"{o}bius
transformations or (c) the Descartes configurations given
by an orbit of the Apollonian group.
The first two of these notions do not generalize to
all dimensions, but version (c) does.

\begin{defi}
{\em 
An {\em Apollonian cluster ensemble} in $n$ dimensions is
defined to be the  orbit of the Apollonian group $\sA_n$ 
 of a given  
ordered, oriented Descartes configuration $\sD_0$
in the parameter space $\sM_{\DD}^n$. }
\end{defi}

This object can
be viewed as a discrete set in the parameter space $\sM_{\DD}^n$
of all (ordered, oriented) Descartes configurations in
dimension $n$. In  dimension
$n=3$ the spheres in an Apollonian
cluster ensemble  yield a sphere packing.
However  in dimensions $n \ge 4$ the spheres in any
such ensemble overlap and no longer correspond 
to a packing, as shown  by Lemma~\ref{le4.1.1}.
Because the Apollonian group consists of rational
matrices, we can ask if there are $n$-dimensional
Apollonian cluster ensembles with rationality
properties, either of their curvatures or of their
full augmented curvature-center coordinate matrices $\bW_{\sD}$.

Similarly one can define
a {\em dual Apollonian cluster ensemble} to be  an orbit of the
dual Apollonian group $\sA_n^{\perp}$ of a single
Descartes configuration,
It is necessarily 
a discrete set, because $\sA_n^{\perp}$ has an integral structure,
so the identity matrix is isolated in the group  $\sA_n^{\perp}$.
We can ask about rationality or integrality properties of
orbits of this group. 

A {\em super-Apollonian cluster
ensemble} is an orbit of the super-Apollonian group  $\sA_n^{S}$. 
We can ask about rationality properties of orbits of this
group, and integrality properties in dimensions $n=2$ and $n=3$.
For $n=2$ these were answered in part II, and for $n=3$ we
show in \S6 that integral $\bW_{\sD}$ do not exist.
For $n \ge 4$ these  orbits are not discrete, but in certain
dimensions 
can have a rational structure inherited from the 
super-Apollonian group, see  \S6.

%
%

\section{Presentation for $n$-Dimensional Apollonian Group}
\setcounter{equation}{0}
In \cite{GLMWY21} we obtained a presentation for the
super-Apollonian group $\sA_{2}^{S}$ in dimension $2$,
which established that it was a hyperbolic Coxeter group.
Theorem~\ref{le4.2.1} gave some non-trivial relations
in the super-Apollonian group $\sA_{n}^{S}$,
all of Coxeter type.  There are additional
relations, at least in dimension 3, see Theorem~\ref{th51}.
It may well be that $\sA_{n}^{S}$ is a hyperbolic Coxeter
group for all $n \ge 3$. However it 
seems a complicated problem to determine a presentation
of  $\sA_{n}^{S}$ in general, and here we
establish a more limited result.
 
We  determine a presentation for the Apollonian group $\sA_n$
in dimensions $n \ge 3,$  given in the next two results.
In all cases it is a hyperbolic Coxeter group.

\begin{theorem}~\label{th51}
For dimension $n=3$ the Apollonian group
$\sA_3 = \langle \bS_1, \bS_2, \bS_3, \bS_4, \bS_5 \rangle$ 
has the relations $\bS_j^2 = \bI$ for $1 \le j \le 5$
and the additional relations
\beql{N501} 
(\bS_j \bS_k)^3 = \bI, \mbox{when}~~ j \ne k.
\eeq
These are a generating set of relations, and $\sA_3$
is a hyperbolic Coxeter group.
\end{theorem}
\paragraph{Proof.}
Recall that
\beql{N539}
\bS_1 = \left[ \begin{array}{rrrrr}
-1 & 1 & 1 & 1 & 1 \\
0 & 1 & 0 & 0 & 0 \\
0 & 0 & 1 & 0 & 0 \\
0 & 0 & 0 & 1 & 0 \\
0 & 0 & 0 & 0 & 1 
\end{array} \right].
\eeq
and the other $\bS_j$ are permutations of the $j$-th and first rows and columns
of $\bS_1$.
It is easy to check that the generators satisfy all the relations 
given above; in what follows we denote this set of relations 
$\sR.$  We note that relations \eqn{N501} can be written in
either of the forms $(\bS_j \bS_k)^2 = \bS_k \bS_j$ or
$\bS_j \bS_k \bS_j = \bS_k \bS_j \bS_k$.

We now consider a general word $U = \bU_1 \bU_2 \cdots \bU_n$
with each $\bU_j = \bS_{i_j}$, where inverses of generators
are eliminated by relations  
 $\bS_i^{-1} = \bS_{i}$. A {\em subword} of $\bU$ is any
word $\bU_{j} \bU_{j+1} \cdots \bU_k$.
We call a
word {\em reduced} if it has the following properties.

(i) {\em $2$-reduced.} 
It contains no subword of form $\bS_j \bS_j$.

(ii) {\em $B$-reduced.}
It contains no subword of form $\bV_1 \bV_2 \cdots \bV_{2m}$
in which $\bV_1 = \bV_3$, $\bV_{2j} = \bV_{2j+3}$ for $1 \le j \le m-2$,
and $\bV_{2m-2}=\bV_{2m}$.

Conditions (i) and (ii) together allow the refinement
of (ii) to assert $\bV_1 = \bV_3 \ne \bV_2$,
$\bV_{2j} = \bV_{2j+3} \ne \bV_{2j+2}$, $\bV_{2m-2}=\bV_{2m}\ne \bV_{2m-1}$.
Also the case $m=2$ rules out words $\bS_j \bS_k\bS_j \bS_k$.
The definition also implies: every subword of a $B$-reduced
word is also $B$-reduced.

We assert that every non-reduced word can be simplified to a reduced
word of a shorter length, using the relations $\sR$.
If (i) is violated, then  using the
relation $\bS_j^2 = I$ we may replace the subword by the empty word, 
shortening
it by $2$.
If (ii) holds then we can replace the subword by
$ \bV_2 \bV_3 \cdots \bV_{2m-1}$, decreasing its length by two. 
This reduction is done by a sequence of replacements of
of the form $\bS_j \bS_k \bS_j$ by  $\bS_k \bS_j \bS_k$ 
applied successively at positions $1, 3, ...,
2m-3$, which ends with
 $\bV_2$ at the beginning and $\bV_{2m-1} \bV_{2m} \bV_{2m}$
at the end. Then the final letters $\bV_{2m} \bV_{2m}$ are deleted
by rule (i), achieving a shortening by $2$. Starting
with any word $\bU$ and applying  this
reduction process, we eventually arrive at a reduced word
or the empty word.

The theorem 
is equivalent to the assertion that each  nonempty  reduced word
is not the  identity. For if the relations $\sR$ did not generate
all relations, there would exist another
relation, necessarily forcing some  nonempty reduced
word to be the identity.

We introduce some (matrix) invariants associated to a word in the
generators. Let $\mathbf{e}_j$ be the $j$-th unit (column) vector,
i.e. the $j$-th column of the $5 \times 5$
identity matrix,  and set
$\mathbf{1}_5 = (1,1,1,1,1)^T$. We define
\beql{N541}
\sigma_j(\bU) := \mathbf{e}_j^T \bU \mathbf{1}_5, 
\eeq
the sum of entries in row $j$ of $\bU$, and 
\beql{542}
\Sigma(\bU) :=\mathbf{1}_5^T  \bU \mathbf{1}_5,
\eeq
the sum of all entries of $\bU$, which we call its {\em mass}. 

We compute the effect on these invariants
of multiplying by a generator.
The matrix $\bS_j \bU$
has all rows but the $j$-th row the same as $\bU$, with
its  $j$-th row equal to the sum of all rows of $\bU$
minus twice its $j$-th row, which yields 
\beql{N543}
\sigma_j (\bS_j \bU) = \Sigma(\bU) - 2 \sigma_j(\bU),
\eeq
and $\sigma_k (\bS_j \bU) = \sigma_k(\bU)$ if $k \ne j$. It also gives
\beql{N544}
\Sigma(\bS_j \bU) = 2 \Sigma(\bU) - 3 \sigma_j(\bU).
\eeq

Now, we define
\beql{N545} 
\delta_j(\bU) := \Sigma( \bS_j \bU) - \Sigma (\bU).
\eeq
This quantity measures the increase in the total mass of a matrix
when we multiply by $\bS_j$.
Two key properties of this measure 
are that if $j \ne k$ then 
\beql{N546a} 
\delta_j(\bS_k \bU) = \delta_j(\bU) + \delta_k(\bU), 
\eeq
while if $j=k$ then 
\beql{N546b} 
\delta_j(\bS_j \bU) = -\delta_j(\bU).
\eeq
To verify these, first observe that 
$$
\delta_j(\bU) = (2\Sigma(\bU) - 3 \sigma_{j}(\bU)) - \Sigma(\bU)
= \Sigma(\bU) - 3 \sigma_{j}(\bU).
$$
Therefore
\begin{eqnarray*}
        \delta_j(\bS_k \bU) &=& \Sigma(\bS_k \bU) - 3\sigma_j(\bS_k \bU) \\
                     &=& 2\Sigma(\bU) - 3\sigma_{k}(\bU) - 3\sigma_j(\bU) \\
                     &=& \Sigma(\bU) - 3\sigma_j(\bU) + \Sigma(\bU) 
- 3\sigma_k(\bU) \\
                     &=& \delta_j(\bU) + \delta_k(\bU),
\end{eqnarray*}
which gives \eqn{N546a}. Similarly, using \eqn{N545},
\begin{eqnarray*}
\delta_j(\bS_j \bU) &=& \Sigma(\bS_j\bS_j \bU) - \Sigma(\bS_j\bU) \\
&=& -(\Sigma(\bS_j\bU)- \Sigma(\bU) ) \\
 &=& -\delta_j(\bU).
\end{eqnarray*}
One consequence of these two properties  is that 
\beql{N547a}
\delta_j (\bS_k \bU) = \delta_k(\bS_j \bU)~~\mbox{if}~~ j \ne k.
\eeq
Another consequence is  
\beql{N548}
\delta_j(\bS_k \bS_j \bU) =
\delta_k(\bS_j \bU) + \delta_j(\bS_j \bU) = 
(\delta_j(\bU) + \delta_k(\bU)) - \delta_j(\bU)
= \delta_{k}(\bU).
\eeq

We assert that all nonempty reduced words $\bU = \bS_j \bU'$ have 
\beql{N549}
\delta_j(\bU') > 0, 
\eeq
so that $\Sigma (\bU) > \Sigma(\bU').$
If this is proved, then since $\bU'$ is also reduced, 
we obtain  by  induction on the length of $\bU$ that
$$ 
\Sigma(\bU) \ge \Sigma(\bU_n) = 7 > \Sigma(I) = 5, 
$$
so that $\bU$ cannot be the identity matrix, and the 
theorem follows. 

We establish \eqn{N549}
 by induction on the length $n$ of 
$\bU= \bU_1 \bU_2 \cdots \bU_n = \bU_1 \bU'$,
where we suppose $\bU_1 = \bS_j$. In the  base case  $n=1$,
we consider that $\bU'=I$,  and we then have
$$
\Sigma(\bS_j) = 7 \ge \Sigma(I) = 5,
$$
completing the base case.
Now suppose $n \ge 2$ and  that the induction
hypothesis holds up to $n-1$.  
We write  $\bU = \bS_j \bS_k \bU''$
with $\bU''$ of length $n-2$, noting that $j \ne k$.
We are to show $\delta_j(\bU') >0$.
Now \eqn{N546a} gives
\beql{N550}
\delta_j(\bU') = \delta_j(\bS_k \bU'') =
\delta_j(\bU'') + \delta_k(\bU''). 
\eeq
The induction hypothesis gives $\delta_k(\bU'')>0$. 
If  $\bS_j \bU''$ is reduced, then 
$\delta_j(\bU'') >  0$ by the induction hypothesis
and we are done. 
So suppose $\bS_j \bU''$ is not reduced. 
Since 
$\bU''$ is reduced, a non-reduced subword in it must
be an initial segment, which either fails to be
$2$-reduced or $B$-reduced.

Suppose first that $\bS_j \bU''$  is not $B$-reduced,
having an initial $B$-word $\bV_1 \cdots \bV_{2m}$ which
begins $\bS_j \bS_l \bS_j$ so that 
$\bS_j \bU'' = \bS_j \bS_l \bS_j \bV$, for some $l \ne j$.
We also have $l \ne k$ since $\bU' = \bS_k \bU''$ is $2$-reduced.
Now we have,  using \eqn{N547a}, 
\beql{N551}
\delta_j(\bU') = \delta_j(\bS_k \bU'') =  \delta_k(\bS_j \bU'').
\eeq
Applying the $B$-reduction
procedure to the initial segment $\bV_1 \cdots \bV_{2m}$ 
simplifies the word  $\bS_j \bU''$ to a word 
$\bS_l \bS_j \bV'$ that is shorter by two letters than  $\bS_j \bU''$,
but is equal to it as a matrix. 
It follows using \eqn{N551} that 
\beql{N552}
\delta_k(\bS_l \bS_j \bV') \equiv \delta_k(\bS_j \bU'') = \delta_j(\bU').
\eeq
We now assert that
 $\bS_l \bS_j \bV'$
itself is reduced. To see this, write
$\bS_l \bS_j \bV' = \bV_2 \bV_3 \cdots \bV_{2m-1} \bV''$
with $\bV''= \bS_r \bS_t \bV'''$. Since $\bU''=\bS_j \bV'$ 
was reduced, any 
further reduction involves a subword $\bX$ that contains $\bV_{2m-2} \bS_r$.
Now $\bS_r \ne \bV_{2m}$ or else $\bU''$ includes
$\bV_{2m} \bS_{r} = \bV_{2m} \bV_{2m}$ contradicting $\bU''$ being reduced.
Next $\bS_r \ne \bV_{2m-1}$, for otherwise  $\bU''$
would contain the subword  $\bV_{2m-2} \bV_{2m-1} \bV_{2m} \bV_{2m-1}$
which is a $B$-word since $\bV_{2m} = \bV_{2m-2}$,
contradicting $\bU''$ being reduced. It follows that  
$\bV_{2m-2} \bV_{m-1} \bS_r
= \bV_{2m} \bV_{2m-1} \bS_r$ consists of three distinct letters
(generators).
Thus any  non-reduced subword $\bX = \bX'\bV_{2m} \bV_{2m-1} \bS_r \bX''$
inside  $\bS_l \bS_j \bV'$ must be  a $B$-word.
But then   $\bU''$ contains the
$B$-word $\bV_{2m} \bV_{2m-1} \bV_{2m} \bS_r \bX''$, 
contradicting $\bU''$ being
reduced. We conclude that  $\bS_l \bS_j \bV'$ is reduced. 
We deduce that  $\bS_k \bS_l \bS_j \bV'$ is reduced, because
$\bS_l \bS_j \bV'$ is reduced and $\bS_k$ differs from both of
its first two letters. The induction hypothesis now applies
to give $\delta_k( \bS_l \bS_j \bV') > 0$, which with \eqn{N552}
gives $\delta_j(\bU')>0$.

The remaining case is that when  $\bS_j \bU''$ is not $2$-reduced. Then
$\bU'' = \bS_j \bV$, and $\bU = \bS_j \bS_k \bS_j \bV$.
We have then, by \eqn{N548}, that
$$
\delta_j(\bU') = \delta_k(\bV).
$$
If $\bS_k\bV$ is reduced, then the induction hypothesis gives
$ \delta_k(\bV)>0$, which gives the desired result. If it
is not reduced, then it is either not $2$-reduced or
$B$-reduced. If it is not $2$-reduced, then we have
$\bV = \bS_k \bV'$, in which case 
the original word $\bU = \bS_j \bS_k \bS_j \bS_k \bV'$
is not $B$-reduced, contradicting the hypothesis. 
If $\bS_k\bV$ is not $B$-reduced, then we have 
$\bV = \bS_l \bS_k \bV'$ with $l \ne k$, and
$\bU = \bS_j \bS_k \bS_j \bS_l \bS_k \bV'$ which shows that $l \ne j$.
Now we have by \eqn{N548} again that 
$$
\delta_k(\bV) = \delta_k (\bS_l \bS_k \bV')
=  \delta_l(\bV') .
$$
If $\bS_k \bV'$ is reduced then the induction hypothesis
gives $ \delta_l(\bV')>0$ which implies that $\delta_j(\bU')>0$
as desired. If it is not $2$-reduced then $\bV' = \bS_k \bV''$
and this contradicts $\bU$ being reduced. There remains
the case where $\bV'$ is not $B$-reduced. Then we get
$\bV' = \bS_m \bS_l \bV''$, with $\delta_k(\bV') = \delta_m(\bV'')$.
We can continue recursively in this way until the entire
word $\bU$ is used up and there are not enough letters to
have any $B$-reduced word. We then  obtain
$$
\delta_j(\bU') = \delta_k(\bV) = \delta_l(\bV) =
\cdots = \delta_n(\bV^{(r)}) > 0. 
$$
This completes the induction step, proving \eqn{N549},
and the theorem follows.
$~~~\bsq$

%
%
\begin{theorem}~\label{th52}
For dimension  $n \ge 4$ the Apollonian group 
$\sA_n = \langle \bS_1, \bS_2, ..., \bS_{n+2} \rangle$ is a 
hyperbolic Coxeter group whose only relations are
$$ 
\bS_j^2 = \bI_{n+2} ~~~\mbox{for}~~~ 1 \le j \le n+2.
$$
\end{theorem}

\paragraph{Proof.} 
Suppose $n \ge 4.$
We must show that 
no nonempty product of the $\bS_j$'s, with distinct
adjacent terms, is equal to the identity matrix $I_{n+2}$.

Let $g_n = \frac{1}{n-1}$ and let $\bA_j =
\frac{1}{2}(\bS_j - \bI_{n+2})$ so that $\bS_j = \bI_{n+2} + 2\bA_j$.  
Then $\bA_j$ has
all entries zero except for the $j$th row, in which all
entries are  $g_n$ except for the
$j$th element, which is $-1$.  Let 
$$
\bA = g_n\bo_{n+2}\bo_{n+2}^T - n \bI_{n+2}
$$ 
be the matrix with $-1$ on the main diagonal and $g_n$
elsewhere, so that its $(i,j)$ element is
\begin{equation*}
a(i,j) = \{ \begin{array}{cl}
 -1, &  \mbox{if}~~i = j, \\
 g_n, &  \mbox{if}~~i \ne j.
\end{array} 
\end{equation*}
Then $\bA_j$ and $\bA$ both have $j$th row $(a(j,1), \ldots, a(j,n+2))$.
We note that $\bA$ is nonsingular.

Suppose that $k \ge 2$, that $1 \le j_1, \cdots, j_k \le n+2$, that no
two consecutive $j_r$'s are equal, and that $S$ is the product $S =
\bS_{j_1} \bS_{j_2} \cdots \bS_{j_k}$.  We assume that $\bS = \bI_{n+2}$ and
derive a contradiction.  We have:
\begin{eqnarray*}
\bS - \bI_{n+2}  &=& \bS_{j_1} \bS_{j_2} \cdots \bS_{j_k} - \bI_{n+2} \\
&=& (\bI_{n+2} + 2\bA_{j_1}) \cdots (\bI_{n+2} + 2\bA_{j_k}) -\bI_{n+2} \\
&=& 2\sum_{1 \le r \le k} \bA_{j_r}
          + 4\sum_{1 \le r < s \le k} \bA_{j_r} \bA_{j_s}
          + 8\sum_{1 \le r < s <t \le k} \bA_{j_r} \bA_{j_s} \bA_{j_t}
          + \cdots
          + 2^k \bA_{j_1} \cdots \bA_{j_k}.
\end{eqnarray*}

The $\bA_j$'s multiply more simply than the $\bS_j$'s, compensating for the
more complicated expression involving them.  A product like $8 \bA_{j_r}
\bA_{j_s} \bA_{j_t}$, which we call an {\em $A$-product}, has one nonzero row,
row $j_r$, which is equal to row $j_t$ of $\bA$ multiplied by $8 a(j_r,j_s)
a(j_s,j_t)$.  The terms of the latter product are $-1$ or $g-n$, depending
on whether or not $j_r = j_s$ and $j_s = j_t$.  We call 
the scalar $8 a(j_r,j_s)a(j_s,j_t)$ 
the {\em $a$-product} corresponding to the $A$-product 
$8\bA_{j_r}\bA_{j_s} \bA_{j_t}$.  In general, 
an $a$-product has the form 
$\pm 2^\alpha g_n^\beta, \alpha > \beta$.

The $j$th row of $\bS-\bI_{n+2}$ gets a contribution from each $A$-product
that starts with $\bA_j$.  It is a linear combination of rows of $\bA$
determined by the last subscript in each contributing $A$-product.
Because this linear combination is zero, each of its coefficients must
be zero, since $\bA$ is full rank.  Therefore, for each $t$ and $u$ in
$\{j_1,...,j_k\}$, the sum of all the $A$-products corresponding to
$A$-products that start with $\bA_t$ and end with $\bA_u$ must be zero.
This sum is a polynomial in the variable $g= g_n$, 
call it $P_{t,u}(g)$.  Because
$A$-products have the form $\pm 2^\alpha g^\beta, \alpha > \beta$,
the coefficient of $g_n^i$ in this polynomial is divisible by $2^{i+1}$, so
we may write 
$$
P_{t,u}(g) = \sum_{i=0}^{k-1} c_{t,u,i} 2^{i+1} g^i
$$ 
for integers $c_{t,u,i}$.

Now consider $P_{j_1,j_k}(g)$. This has degree $k-1$.  In fact, the
only $A$-product that can contribute an $a$-product of degree $k-1$ is
$2^k \bA_{j_1} \cdots \bA_{j_k}$, and since successive $\bA_{j_r}$'s are
distinct, the corresponding $a$-product is, $2^k g_n^{k-1}$, i.e.,
$c_{j_1,j_k,k-1} = 1$.  Therefore, by the Rational Root Theorem, since
%
%
$g_n = \frac{1}{n-1}$ is a root of $P_{j_1,j_k}(g)$, we must 
either have $n-1$ dividing $2^k$,
or $n = 2^m + 1$ for some $m > 1$, in which case, $g_n = \frac{1}{2^m}$.
Writing $c_i = c_{j_1,j_k,i}$ we then have
$$
P_{j_1,j_k}(g_n) = \sum_{i=0}^{k-1} c_i 2^{i+1} g_n^i = 0.
$$
On multiplying by $2^{m(k-1)-k}$ this 
yields the equation
$$
\sum_{i=0}^{k-1} c_i 2^{(m-1)(k-1-i)} = 0, 
$$
with integer entries $c_i$.
We now check that this equation is impossible $(\bmod~2)$.
All terms except the last are even since $(m-1)(k-1-i) > 0$ when
$0 \le i \le k-2$. The 
$(k-1)$-st  term in the last sum is $c_{k-1} 2^0 = 1$, so we
get $0 \equiv 1~(\bmod~2)$, 
the desired  contradiction.  We conclude that  $\bS \ne \bI_{n+2}$, 
and the theorem follows. $~~~\bsq$

%
%
%
\section{$S$-Integral Apollonian Cluster Ensembles (Dimension $n$)}
\setcounter{equation}{0}

We study integrality and rationality properties for
Apollonian cluster  ensembles. Given a finite set of primes $S$,
we  say that 
a rational number is {\em $S$-integral}
if its denominator is divisible only by powers of primes in $S$. 
By convention we let $S=1$ denote the case when there are no primes in $S$.

The Apollonian group consists of integer matrices in
dimensions $2$ and $3$, and retains an $S$-integral structure
in all dimensions, for suitable $S$. In each dimension $n$
we  consider the questions:

(1) Does there exist some $S$ 
and an Apollonian cluster ensemble
all of whose Descartes configurations consist of spheres 
having  $S$-integral
curvatures (``$S$-integral ensemble'')?

(2) Does there exist some $S$ 
and an Apollonian cluster ensemble
all of whose Descartes configurations have
augmented curvature-center coordinate matrices $\bW_{\sD}$ 
$S$-integral (``super $S$-integral ensemble'')?

We show that $S$-integral ensembles exist in all
dimensions, if $S$ is chosen properly. However
we show that super $S$-integral  ensembles  can exist only
in dimensions $n = 2m^2$ or $n = (2m+1)^2$,
for integer $m$, again with $S$ chosen appropriately.

%
%
%

\subsection{$S$-integral Apollonian Cluster Ensembles}

We say that an 
 Apollonian cluster ensemble is {\em $S$-integral} if
the curvature of every sphere in the ensemble is $S$-integral.

\begin{theorem}~\label{th61} 
In each dimension $n \geq 2$ there
exists an  $S$-integral Apollonian cluster ensemble in which
$S$ is specified as:

(1) $S$ is  the set of primes dividing $n - 1$ if $n$ is
even.

(2) $S$ is the set of primes dividing  $\frac{n-1}{2}$ if $n$ is odd.
\end{theorem}

\pf
It suffices to show that the Descartes equation
\beql{510a} 
 \bb^T \bQ_{D,n} \bb = 0
\eeq
has a non-zero $S$-integral solution $\bb$ for each $n \geq 2.$ 
There  is such a configuration $\sD$ which is not only $S$-integral,
but integral, 
with curvatures $(0, 0, 1, 1,..., 1).$
It consists of two parallel hyperplanes separated by distance $2$
together with $n$ unit spheres whose centers comprise
the vertices  an $(n - 1)$-dimensional
simplex in a hyperplane parallel to
the two hyperplanes in the configuration, and lying
midway between them. 

The other Descartes configurations in the Apollonian cluster ensemble
and super-Apollonian cluster ensemble generated by this configuration
are $S$-integral, where $S$ is the set of  primes
dividing the denominator of $\frac{2}{n - 1},$ since they
have associated matrices  $\bU \bW_{\sD}$ for some $\bU$ in the Apollonian
group. ~~~$\bsq$ \\

In dimension $2$ we can take $S=1$, as
we saw in \cite{GLMWY22}. In that case an
Apollonian cluster
ensemble consists of the Descartes configurations in 
an Apollonian circle packing. 
In \cite{GLMWY2} we studied in various 
number-theoretic questions
related to the integer curvatures  appearing in 
integer Apollonian circle  packings. 

%
%
%

\subsection{Super $S$-Integral Apollonian Cluster Ensembles}

We say that
a Descartes configuration $\sD$  is {\em super $S$-integral}
if its augmented curvature-center coordinate matrix $\bW_{\sD}$ is
$S$-integral. Similarly, an 
Apollonian cluster ensemble is {\em super $S$-integral} if
every Descartes configuration $\sD$ in the packing
has  an  $S$-integral augmented curvature-center matrix. 
The next lemma reduces the  question of super $S$-integrality of
an Apollonian cluster ensemble to that of a single
Descartes configuration.

%
%
%

\begin{theorem}~\label{le62}
If a single
Descartes configuration is super 
$S$-integral, then the Apollonian cluster ensemble it generates is
super  
$S^\prime$-integral,
where $S^\prime$ consists of $S$ together with all primes
dividing the denominator
of $\frac{2}{n-1}$.
\end{theorem}
\pf
This follows from the fact that the Apollonian group 
consists of rational matrices whose entries have denominators
that are divisible only by primes dividing the
denominator of  $\frac{2}{n-1}$, when put in
lowest terms.  $~~~\bsq$ \\

It seems to be a difficult problem to determine
for specific $S$ for which dimensions there might
exist a super $S$-integral Descartes configuration.
We consider the weaker question of whether
in a given dimension there exists a  super $S$-integral Descartes
configurations for some $S$. This is the same as the existence of
Descartes configurations $\sD$ having a rational ACC-matrix
$\bW_{\sD},$ and we call such Descartes configurations
{\em super-rational.}

According to the augmented Euclidean Descartes Theorem~\ref{th31},  
super-rational Descartes configurations occur exactly
in those dimensions $n$ in which there exists an
invertible rational matrix $\bW$ such that 
\beql{505b}
\bW^T  \bQ_{D,n} \bW = \bQ_{W,n} :=  
  \left[ \begin{array}{ccl}
      0     &  -4    &   0    \\
     -4     &   0    &   0   \\
      0     &   0    &   2I_n    \end{array} \right],
\eeq
that is,  the quadratic form $Q_{D,n}$ is rationally
equivalent to the form $Q_{W,n}$ 
We use this fact to determine in which dimensions such 
configurations exist.

\begin{theorem}~\label{th63}
A necessary and sufficient condition on the dimension $n$
for a super $S$-integral  Descartes
configuration to exist for some $S$ is that 
$n = 2k^2$ or $(2k - 1)^2$
for some positive
integer $k$.
\end{theorem}

To establish this result, we proceed in a series of  lemmas.

\begin{lemma}~\label{le64}
Given a Descartes configuration $\sD$ in $\rr^n$ its associated
augmented matrix $\bW_{\sD}$ has 
\beql{501}
\det (\bW_{\sD})^2 = n 2^{n + 3}.
\eeq
\end{lemma}

\pf
This follows from taking determinants in \eqn{307},
since the right side has determinant $- 2^{n+4}$
while the left side has determinant 
$\det (\bW_{\sD})^2 \det(\bQ_{\sD, n})$ and 
\beql{502}
\det (\bQ_{\sD,n}) = -\frac{2}{n}.
\eeq
To verify this last statement, we apply the following
row operations to the matrix $\bQ_{n}$. Add rows $2$ through $n+2$
to the first row, to get a new first row that has all entries
$-\frac{2}{n}$. Then add this row multiplied by $-\frac{1}{2}$
to each of the other rows. Aside from the first row, the
first column is zero, and the lower right $(n+1) \times (n+1)$
matrix is the identity. But this matrix obviously has
determinant $-\frac{2}{n}.$
~~~$\bsq$

\begin{lemma}~\label{le65}
If a super-rational  Descartes
configuration exists in dimension $n$, then necessarily
$n = 2k^2$ or $(2k - 1)^2$
for some positive
integer $k$.
\end{lemma}

\pf
A necessary condition  
for the existence of a Descartes configuration
$\sD$ whose augmented  matrix  
$\bW_{\sD}$ has rational entries
is that
$\det (\bW_{\sD})$ be rational. 
This  requires
that $ n2^{n + 3}$ be the square of a rational number.
By Lemma~\ref{le64}, this holds for even $n$ if and only if
$n$ is twice a square, 
and for odd $n$ if and only if $n$ is an (odd) square.~~~$\bsq$

To prove the sufficiency of 
this condition, we use the theory
of equivalence of rational quadratic forms, cf. Cassels~\cite{Cas78}
or Conway~\cite{Con97}.  We write
$\bQ \simeq_{\qq} \bQ^\prime$ to mean that 
the (rational) quadratic form $\bQ$ is rationally equivalent to 
$\bQ^\prime$. To apply
the decision procedure, we first  diagonalize 
$\bQ_n$ over the rationals, which we do for all $n \geq 2.$
%
%

\begin{lemma}~\label{le66}
For each $n \geq 2$, the Descartes quadratic form 
$\bQ_{D,n} = \bI_{n + 2} -\frac{1}{n} \bo_{n+2} \bo_{n+2}^T$ has
\beql{505c}
\bQ_{D,n} ~\simeq_{\qq}~ \mbox{diag} (\frac{n-1}{n}, \frac{n-2}{n-1}, 
\cdots, \frac{2}{3}, 2, 2, 2, -2).
\eeq
\end{lemma}

\noindent{\bf Proof.}
We diagonalize the quadratic form as in Conway \cite[pp. 92--94]{Con97}.
Set
$$M^{(n+2)} := \bQ_{D,n} =
(x_0 + y_0 ) \bI_{n+2} - y_0 \bo_{n+2} \bo_{n+2}^T \,,
$$
where $x_0 = \frac{n-1}{n}$,
$y_0 = \frac{1}{n}$.
At the $j$-th stage of reduction we will have
$$ \bQ_{D,n} \simeq_\qq {\rm diag} (d_1, d_2, \ldots, d_j, 
M^{(n+2-j)} ) \,,
$$
where
\beql{VP428}
\bM^{(n+2-j)} = (x_j + y_j ) \bI_{n+2-j} - y_j \bo_{n+2-j} \bo_{n+2-j}^T
\eeq
for certain $x_j$, $y_j$.
The reduction step is
\beql{VP429}
( \bW^{(j)} )^T \bM^{(n+2-j)} \bW^{(j)} =
{\rm diag} (d_{j+1} , \bM^{(n+1-j)} ) \,.
\eeq
To specify $\bW^{(j)}$ we first let
 $\bW_m ( \alpha )$ be the $m \times m$ real matrix
$$\bW_m ( \alpha ) =
\left[ \begin{array}{cc}
1 & \alpha \cdots \alpha \\
~ & ~ \\
{\bf 0} & \bI_{m-1}
\end{array}
\right] \,, 
$$
and we set
\beql{VP430}
\bW^{(j)} := \bW_{m+2-j} \left( \frac{y_j}{x_j} \right).
\eeq
Substituting this in \eqn{VP429}, its left side yields a matrix 
with the form of the right side with
$$d_{j+1} = x_j \,,$$
and with  $x_{j+1}$, $y_{j+1}$ given by the recursion
\begin{eqnarray}\label{VP431}
y_{j+1} & = & y_j + \frac{y_j^2}{x_j} \,, \\
\label{VP432}
x_{j+1} + y_{j+1} & = & x_j - \frac{y_j^2}{x_j} ~.
\end{eqnarray}
Solving this recursion, by induction on $j$, one obtains
$$
\begin{array}{rllcl}
x_j & = & \displaystyle\frac{n-j-1}{n-j} \,, &~~~~ & 0 \le j \le n-2, \\ [+.2in]
y_j & = & \displaystyle\frac{1}{n-j} \,, && 0 \le j \le n-2 \,.
\end{array}
$$
This yields the diagonal elements
\beql{VP433}
d_j = \frac{n-j-1}{n-j}, \quad
1 \le j \le n-3 \,,
\eeq
with
$$ \bQ_{D,n} \simeq_\qq~
{\rm diag} ( \frac{n-1}{n} , \ldots, \frac{2}{3} , d_2, \bM^{(4)} ) \,.
$$
We find
$d_2 = x_3 = \frac{2}{3}$ and
$$
\bM^{(4)}
= (x_{n-2} + y_{n-2} ) \bI_4 -
y_{n-2} \bo_4 \bo_4^T = \frac{1}{2}
\left[ \begin{array}{rrrr}
1 & -1 & -1 & -1 \\
-1 & 1 & -1 & -1 \\
-1 & -1 & 1 & -1 \\
-1 & -1 & -1 & 1
\end{array}
\right] = \bQ_{D,2} \,.
$$
For the final step in the reduction we use
\beql{VP434}
\bN^T \left( \bQ_{D,2} \right) \bN =
{\rm diag} (2,2,2,-2), 
\eeq
with
$$
\bN = \left[ \begin{array}{rrrr}
1 & -1 & -1 & 1 \\
-1 & 1 & -1 & 1 \\
-1 & -1 & 1 & 1 \\
1 & 1 & 1 & 1
\end{array}
\right] \,.
$$
This completes the reduction.~~~$\bsq$ \\

%
%

\noindent{\bf Proof of Theorem~\ref{th63}.}
The necessity 
of $n = 2k^2$ or $(2k - 1)^2$ was proved in Lemma~\ref{le65}.
The sufficiency is equivalent to proving that if $n = 2k^2$
and $n = (2k - 1)^2$ then 
$$
\bQ_{D,n} \simeq_{\qq} \bQ_{W,n} :=   
\left[ \begin{array}{rcr}
           0 &  -4      &    0   \\
          -4 &   0      &    0   \\
           0 &   0      &   2\bI_n    \end{array} \right]. 
$$
We begin by noting the rational equivalence 
\beql{VP427}
\bQ_{W,n} \simeq_{\qq}  \mbox{diag}(-2, 2, \cdots, 2, 2)
= {\rm diag} (-2, 2I_{n+1})
\eeq
via the matrix
$$
\bW_0 = \frac{1}{2} \left[ \begin{array}{ccl}
            1  &    1      &   0   \\
            1  &   -1      &   0   \\
            0  &    0      &   2\bI_n    \end{array} \right].
$$
By permuting variables we have 
$\bQ_{W,n} \simeq_{\qq} \mbox{diag}(2, 2, \cdots, 2, - 2).$
Thus the theorem is equivalent to showing that
$\bQ_{D,n}$ is rationally equivalent to $\mbox{diag}(2, 2, 2,..., -2)$.
Lemma \ref{le66} gives
\begin{eqnarray}\label{VP435}
\bQ_{D,n} & \simeq_\qq & 
( \frac{n-1}{n} , \frac{n-2}{n-1}, \ldots,
\frac{3}{2} , 2,2,2,-2) \,, \nonumber \\
& \simeq_\qq &
(n(n-1), (n-1) (n-2) , \ldots, 3 \cdot 2, 2,2,2, -2) \,,
\end{eqnarray}
using at the last step a conjugacy by
$\bW = {\rm diag} (n, n-1, \ldots, 2,1,1,1,1)$.

The Hasse-Minkowski theorem says that two rational quadratic
forms of the same dimension are equivalent if and only
they have the same signature, 
the ratio of their  determinants is a nonzero square, and
they are $p$-adically equivalent for all primes $p$,
 cf. Conway~\cite[p. 96ff]{Con97}. Lemma~\ref{le66} shows
that the signatures of $\bQ_{D,n}$ and $\mbox{diag}(2, 2, 2,..., 2, - 2)$
agree, and the hypothesis $n = 2k^2$ or $n = (2k - 1)^2$ is
exactly the condition that the ratio of their determinants 
is a square of a rational number, and it remains to check the
$p$-adic invariants.

The $p$-adic invariants $\sigma_p (\bQ )$ are defined $(\bmod~8)$,
and for a diagonal form \\
$\bQ = {\rm diag} (d_1, d_2, \ldots, d_n )$, one has
\beql{VP436}
\sigma_p ( \bQ ) \equiv \sum_{j=1}^n \sigma_p (d_j) ~~(\bmod~8) \,.
\eeq
We recall formulas for $\sigma_p (d)$ when $d \in \zz$, cf.
Conway~\cite[pp. 94--96]{Con97}.
Write $d= b p^l$ with $(b,p) =1$.
For $p \ge 3$, and an even power $l=2j$,
\beql{VP437}
\sigma_p (d) \equiv  p^{2j} \equiv 1 ~~ (\bmod~8) \,,
\eeq
while for an odd power $l=2j+1$,
\beql{VP438}
\sigma_p (d) \equiv \left\{
\begin{array}{llll}
p & (\bmod~8) & \mbox{if} & \left( \frac{b}{p} \right) =1 \,, \\ [+.1in]
p+4 & (\bmod~8) & \mbox{if} & \left( \frac{b}{p} \right) = -1 \,.
\end{array}
\right.
\eeq
If $p=2$ then for an even power $l= 2j$,
\beql{VP439}
\sigma_2 (d) \equiv b ~~ (\bmod~8) \,,
\eeq
while for an odd power $l= 2j+1$,
\beql{VP440}
\sigma_2 (d) \equiv \left\{
\begin{array}{llll}
b & \mbox{if} & b \equiv \pm 1 & (\bmod~8) \,, \\ [+.1in]
b+4 & \mbox{if} & b \equiv \pm 3 & (\bmod~8) \,.
\end{array}
\right.
\eeq
Now \eqn{VP435} gives
$$
\sigma_p ( \bQ_{D,n} ) \equiv
\sum_{j=0}^{n-3}
\sigma_p ((n-j) (n-j-1)) + 3 \sigma_p (2) + \sigma_p (-2) ~~( \bmod~8), 
$$
while \eqn{VP427} gives
$$
\sigma_p ( \bQ_{W,n} ) \equiv
\sum_{j=0}^{n-3} \sigma_p (2) + 3\sigma_p (2) + \sigma_p (-2) ~~ (\bmod~8) \,.
$$
To show equality of these, it suffices to show that for all $p$,
\beql{VP441}
\sum_{j=0}^{n-3} \sigma_p (2) \equiv \sum_{j=0}^{n-3} \sigma_p ((n-j)(n-j-1)) 
~~(\bmod~8)
\eeq
holds whenever $n= 2k^2$ or $n= (2k-1)^2$.

Consider first the case that $p \ge 3$ is odd.
Then each $\sigma_p (2) =1$, so
\beql{VP442}
\sum_{j=0}^{n-3} \sigma_p (2) \equiv n-2 ~ (\bmod~8) \,.
\eeq
Now if $p$ does not divide $ (n-j ) (n-j-1)$ then
$\sigma_p ((n-j) (n-j-1)) =1$.
The terms divisible by $p$ occur in blocks of two consecutive terms, 
and we claim that if $p$ divides $j$ then
\beql{VP443}
\sigma_p ((j + 1)j ) + \sigma_p (j(j-1)) \equiv 2 ~ (\bmod~8).
\eeq
Suppose $j = bp^l,$ with where $(b,p) = 1$ and $l \geq 1.$
If $l$ is even, both terms on the left side of \eqn{VP443}
are 1 $(\bmod~8)$ by \eqn{VP437}, while if $l$
is odd, then if $p \equiv 1$ $(\bmod~4)$,
then $\left( \frac{-1}{p}\right) =1$, so the two terms both have values $p$
(resp. $p+4$) according as  $\left( \frac{b}{p}\right) =1$
(resp. -1),   and their sum is $2p \equiv 2$ $(\bmod~8)$.
If $p \equiv 3$ $(\bmod~4)$, then
$\left( \frac{-1}{p} \right) =-1$, so exactly one of
$\left( \frac{\pm b}{p} \right)$ takes the value $-1$, 
and the two terms add up to $2p + 4 \equiv 2$ $(\bmod~8)$.
Thus \eqn{VP443} follows.
Thus adding up the right side of \eqn{VP441} and 
grouping terms divisible by $p$ in consecutive pairs gives
\beql{VP444}
\sum_{j=0}^{n-3}
\sigma_p ((n-j) (n-j-1)) \equiv
\sum_{j=0}^{n-3} 1 \equiv n-2 ~ (\bmod~ 8) \,.
\eeq
There remains an exceptional case where $p | n$, in which case $n(n-1)$
is divisible by $p$ and is an un-paired term.
Since $n= 2k^2$ or $(2k-1)^2$, thus $p^l | n$ with $l$ even,
hence $\sigma_p (n(n-1)) =1$ in this case, and \eqn{VP444}
holds.
This establishes \eqn{VP441} for $p \ge 3$.

Now consider the case $p=2$.
Certainly $\sigma_2(2) =1$ so \eqn{VP442} holds.
We claim that
\beql{VP445}
\sigma_2 ((2j+1) 2j) + \sigma_2 (2j(2j-1)) \equiv 0 ~ (\bmod~8) \,.
\eeq
Write $2j = 2^l b$ with $b$ odd, and by checking all possible cases using
\eqn{VP439} and \eqn{VP440}, one verifies \eqn{VP445}.
Suppose $n=2k^2$.
Then in the right side of \eqn{VP441} all terms pair except the first and last,
and \eqn{VP445} yields
\begin{eqnarray*}
\lefteqn{\sum_{j=0}^{n-3}
\sigma_2 ((n-j(n-j-1)) \equiv \sigma_2 (n(n-1)) + \sigma_2 (3 \cdot 2)} \\
&& \quad = \left\{ \begin{array}{rll}
-1+ -1 & \mbox{if} & k \equiv 0 ~ (\bmod~2) \,  \\[+.1in]
1 + -1 & \mbox{if} & k \equiv 1 ~ (\bmod~2)
\end{array}
\right. \\
&& \quad = n-2 ~ (\bmod~8) \,,
\end{eqnarray*}
so \eqn{VP441} holds.
If $n = (2k-1)^2 \equiv 1$ $(\bmod~8)$ then all term pair except 
the last term, and \eqn{VP445} yields
$$\sum_{j=0}^{n-3} \sigma_2 ((n-j)(n-j-1)) = \sigma_2 (3 \cdot 2) 
\equiv -1 ~(\bmod~8) \,,
$$
so \eqn{VP441} holds in this case.~~~$\bsq$ \\

Theorem~\ref{th63} establishes the existence
of super-rational Descartes configurations in the given
dimensions, but does not give a bound for the 
denominators of the rationals appearing in these
configurations. It could be that in certain dimensions 
$n = 2k^2$ and $(2k +1)^2$ there exist strongly
integral Descartes configurations, i.e. ones with $S = 1.$
However, even  
if such configurations exist for some $n > 2$, then the Apollonian cluster
ensemble containing them would not inherit the super-integrality
property, but only super $S^{'}$-integrality as in Lemma~\ref{le62}.
We leave this as an open problem.

%
%
%

\section{$n$-Dimensional Duality Operator}
\setcounter{equation}{0}

In the two-dimensional case we showed
that there exists
a duality operator $\bD \in Aut(Q_{D})$ which took
each Descartes configuration $\sD$ to a new Descartes
configuration $\sD'$ that consists of four 
circles  orthogonal to the original circles, each passing
through three of the  intersection points of the
circles of $\sD$.  
We showed that $\bD$  was contained
in the normalizer of the super-Apollonian group
in two dimensions.

The geometric duality operation based on orthogonal spheres 
generalizes to higher dimensions as follows. 
Given $n+1$ mutually tangent spheres in $n$ dimensions, 
there is a unique
sphere through their points of tangency, and this sphere
is orthogonal to each 
of the given $n + 1$ spheres, as given in the following 
(known) result.
%
%
%
%

\begin{theorem}~\label{th71}
Given  $n+1$ mutually tangent $(n - 1)$-spheres 
$\{ C_i : 1 \le j \le n+1 \}$ in 
$\rr^n$ having disjoint interiors,
there  exists a unique $(n-1)$-sphere $C^\perp$ passing through the 
$\frac{n (n +1)}{2}$ tangency points of these spheres.  
At each
such tangency point the normal to the sphere $C^\perp$ is perpendicular
to the normals of the two spheres $C_i$ and $C_j$ tangent there.
\end{theorem}

\pf
The assumption of disjoint interiors (we allow interior to be
defined as ``exterior'' for one sphere)
is equivalent to all
$\frac{n (n +1)}{2}$ tangency points of the spheres being distinct.
For dimension $n = 2$ there is a unique circle
through any three distinct points. 
However for  $n \geq 3$ the
conditions are over-determined, since $n + 1$ distinct points 
already determine a unique $(n - 1)$-sphere, and the main issue
is existence.  

Both assertions of the theorem are invariant under M\"{o}bius
transformations (which preserve angles), 
and there exists a M\"obius transformation
taking a set of $n+1$ mutually tangent $(n - 1)$-spheres in 
$\rr^n$ having disjoint interiors to any other such set, cf.
 Wilker~\cite[Theorem 3]{Wi81}. Thus it suffices to prove
the result for a single such configuration, and we consider
the configuration of $n+1$ mutually
touching spheres of equal radii whose centers are at the vertices of
a regular $n$-simplex, and  tangency points of the spheres
are the midpoints of its edges.
The first assertion of the theorem  holds in this case
because there is an $(n - 1)$- sphere whose center is at the
center of gravity of this simplex, which passes through
the midpoints of every edge of the simplex. Indeed the
isometries preserving an $n$-simplex fix the center of
gravity and act transitively 
on the edges. Note that for $n =2$ the simplex is an
equilateral triangle and $C^\perp$ is the inscribed circle;
however for $n \geq 3$ the sphere $C^\perp$  is neither inscribed nor
circumscribed about this simplex.

For the second assertion of
the theorem, in  this configuration the sphere $C^\perp$ has each edge of
the $n$-simplex lying in a tangent plane to the sphere; so 
the normal to $C^\perp$ at the midpoint of an edge is perpendicular
to that edge.
Two spheres $C_i$ and $C_j$ intersect at the midpoint of an
edge, and the normal to their tangent planes points along
this edge; thus this normal is perpendicular to the normal to $C^\perp$
there. ~~~$\bsq$ \\

The second assertion in Theorem~\ref{th71} explains why
the sphere $C^\perp$ is termed `` orthogonal.'' 
Thus, given a Descartes configuration of $n+2$ spheres
$C_i$, we get a system of $n+2$
``orthogonal'' spheres
$$\sD^\perp := \{C_{1}^\perp, \ldots, C_{n+2}^\perp \},$$
 where  $C_i^\perp$
is associated to the $n + 1$ spheres
obtained by deleting $C_i$.
When $n = 2$ the new spheres are mutually tangent 
and give a new Descartes
configuration; this gives the ``duality'' operation $\bD$.
For $n \geq 3$, however, the spheres are not mutually tangent.
In fact for all $n$ their curvatures satisfy a relation
similar in form to the  original (two-dimensional)
Descartes relation, namely
\beql{eqN61a}
 \sum_{i=1}^{n+2} \frac{1}{r_i^2} = 
\frac{1}{2}( \sum_{i=1}^{n+2} \frac{1}{r_i} )^2,
\eeq
and not the Soddy-Gossett relation 
\beql{eqN61b}
 \sum_{i=1}^{n+2} \frac{1}{r_i^2} = 
\frac{1}{n}( \sum_{i=1}^{n+2} \frac{1}{r_i} )^2,
\eeq
satisfied by Descartes configurations in $n$-dimensions,
cf. \cite[Theorem 1.2]{LMW02}. 
(We omit a proof of \eqn{eqN61a}.)
In particular, for $n \ge 3$ given a Descartes configuration
$\sD$, the set $\sD^\perp := \{C_{1}^\perp, \ldots, C_{n+2}^\perp \}$ of
orthogonal spheres is {\em not}
a Descartes configuration.

The question arises, are these $n+2$ ``orthogonal''
spheres in any special relation to one another? 
We answer this in terms of an inversive invariant of two arbitrary
(not necessarily tangent) oriented spheres.

%
%

\begin{defi}~\label{de71}
{\em
(i) The
{\em separation}  between two oriented  spheres $C_1$ and $C_2$
with finite radii $r_1$ and $r_2$, and with centers distance $d$ apart, is
\beql{N827a}
\Delta(C_1, C_2) := \frac{d^2 - r_1^2 - r_2^2 }{2r_1r_2}, 
\eeq
provided  both spheres are inwardly oriented or outwardly oriented,
and is otherwise the negative of the right side of this formula.

(ii) The {\em separation} of an oriented
sphere  $C_1$ of finite radius $r_1$ and an oriented  hyperplane
$C_2$ is 
\beql{N827b}
\Delta(C_1, C_2) :=\frac{d}{r_1}, 
\eeq
where $d$ is the (signed) distance
from the center $\bx_1$ of $C_1$ to $C_2$,  measured so that $d \geq 0$
if $\bx_1$ is not in the interior of $C_2$ and $C_1$ is inwardly
oriented, or if $\bx_1$ is in the interior of $C_2$ and 
$C_1$ is outwardly oriented, and $d < 0$ otherwise.

(iii) The {\em separation} between two oriented hyperplanes 
$C_1$ and $C_2$ is
\beql{N827c}
\Delta(C_1, C_2) := - \cos \theta, 
\eeq
where 
$\theta$ is 
the dihedral angle between the designated normals at a point of
intersection. 
}

\end{defi} 

The separation of two spheres is an inversive
invariant, hence an M\"{o}bius invariant; that is,
\beql{827d}
\Delta(\fg(C_1), \fg(C_2)) = \Delta(C_1, C_2),
\eeq
holds for any M\"{o}bius transformation $\fg$.
This concept appears in Boyd~\cite{Bo73},
who introduced the  term {\em separation} for it,
but the concept~\footnote{
The idea of considering such an inversive invariant traces 
back to  work of Clifford \cite{Cl1882} in 1868  and of
Darboux \cite{Da1872} in 1872. 
However, neither Clifford's nor Darboux'  definition 
was  precisely  $\Delta(C_1, C_2)$. 
Clifford defines the  {\em power of two spheres} 
to be the square distance of their centers less the sum of the squares of 
their radii, i.e., $d^2-r_1^2-r_2^2$, and
Darboux also uses the same quantity, \cite[p.350]{Da1872}.
} 
was used earlier
by Mauldon~\cite{Ma62} in 1962, who used the term {\em inclination} to
mean the negative of $\Delta(C_1, C_2),$ and showed it was
an inversive invariant. The absolute value of $\Delta(C_1, C_2)$ was
studied in Beardon \cite[pp.29]{Be83} under the name of {\em inversive product}
of two spheres. 

The separation $\Delta(C_1, C_2)$ of two spheres can be expressed in
terms of their augmented curvature-center coordinates as
\begin{eqnarray}~\label{N501bb}
\Delta(C_1, C_2) & = & \frac{1}{2}\bw(C_1)^T \bK_{n} 
\bw(C_2)      \nonumber \\
& = &-\frac{1}{2}(\bar{b}(C_1) b(C_2) + b(C_1)\bar{b}(C_2) ) +
b(C_1)b(C_2)\sum_{j = 1}^n  x_j(C_1)x_j(C_2),
\end{eqnarray}
where $\bK_n(1)$ is given by
\beq
\bK_n = \left[ \begin{array}{ccl}
           0 &    -1      &   0   \\
          -1 &     0      &   0   \\
           0 &     0      &   2I_n   \end{array} \right].
\eeq
This formula can be proved by a simple algebraic
calculation. Using it,
one can check that  for two tangent  spheres $C_1$ and $C_2$, 
$\Delta(C_1, C_2) = 1$, if (1) $C_1$ and $C_2$   are
externally tangent, and both are inwardly oriented or outwardly 
oriented,  or (2) $C_1$ and $C_2$  are internally tangent and one is
inwardly oriented, the other is outwardly oriented. 
In all other cases two tangent spheres have $\Delta(C_1, C_2)= -1$,   
and orthogonal spheres are those with $\Delta(C_1, C_2) =0.$

 From Theorem~\ref{th71} one obtains
\beql{N502aa}
\Delta( C^\perp, C_j) = 0 \qquad\mbox{for}\qquad  1 \le j \le n+1,
\eeq
and these relations determine $C^\perp$ up to orientation.
It can also be shown that if a set of tangent spheres 
$\{C_1,...,C_{n+1}\}$ have oriented
curvatures $\bb_{n+1} = (b_1,..., b_{n+1})$, 
and centers $\bx_j$, then for either orientation
the  orthogonal sphere $C^\perp$ has 
oriented curvature
$q$ satisfying
\beql{N502ab}
q^2  ~ = ~ \frac{1}{2}
\left(
\frac{1}{n-1}( \sum_{j=1}^{n+1} b_j)^2  - \sum_{j=1}^{n+1} b_j^2 \right),
\eeq
and (oriented) center $\bx$ satisfying
\beql{N502ac}
q \bx = -  \bb_{n+1}(\frac{1}{2} \bQ_{D, n-1})\bC,
\eeq
in which $\bC$ is an $(n+1) \times n$ matrix whose $j$-th row is $b_j\bx_j$.

An  oriented Descartes configuration in $\rr^n$ is
characterized in terms of separation as
a set of $n+2$ oriented spheres each pair of which has
$\Delta(C_i, C_j) = 1$, when $i \neq j.$  Thus such a configuration
has the following property.

%
%
\begin{defi}~\label{de51}
{\em
A  collection of
oriented spheres is {\em equiseparated} if all values $\Delta(C_j, C_k)$
with $j \neq k$ are equal.
}
\end{defi}

The equiseparation
property can also be viewed as an {\em equiangularity} property,
because for two oriented circles that intersect or touch one has
\beql{N504a}
\Delta(C_1, C_2) = - \cos \theta,
\eeq
where $\theta$ is the angle between oriented normals at a
point of intersection of the two circles.
We now show the duality operation preserves equiseparability
in all dimensions. 

\begin{theorem}[Equiseparation Theorem]\label{th72}
Given an oriented  Descartes configuration 
$\sD = (C_1, C_2, ..., C_{n+2})$ in $\rr^n$, if the
dual spheres are properly oriented then the 
(oriented)
dual configuration $(C_1^\perp, C_2^\perp, ..., C_{n+2}^\perp)$
is equiseparated, with 
\beql{N703}
\Delta(C_j^\perp, C_k^\perp) =  \frac{1}{n - 1} 
\qquad\mbox{if}\qquad j \neq k.
\eeq 
\end{theorem}

\pf
In this result, the
orientation assigned to the dual spheres in the
theorem depends on all
$n+2$ spheres in the Descartes configuration, and
the orientation of $C_j^\perp$ cannot be consistently assigned from
the $n+1$ oriented spheres $\{ C_i : i \neq j\}$ alone.
If all $n+2$ spheres $C_j$ are inwardly oriented, then
$n+1$ of the spheres $C_j^\perp$ will be inwardly oriented
and one outwardly oriented, the last being the one of
largest radius. If all but one of the $n+2$ spheres are
inwardly oriented, and one outwardly oriented,
then all $n+2$ spheres $C_j^\perp$
will be inwardly oriented.

Since the result is invariant under inversion, it suffices
to prove it for a single Descartes configuration.
We consider the special 
oriented Descartes configuration where the curvatures are
(0,0,1,1,\ldots, 1).   Here we have two parallel planes, which we take as
$x_1 = \pm 1$, and $n$ unit spheres, all with centers on the plane
$x_1 = 0$.  Their centers form a regular simplex in this plane. 
We may take one of these centers at $(0, \xi ,0,0,\ldots)$ where 
$\xi^2 = 2(n-1)/n$.  Consider the ``orthogonal'' spheres that pass through 
the point $T = (1, \xi ,0,0,\ldots,  0)$.  
There are $n$ such spheres, and all but one of 
them is a hyperplane containing 
the points $T$ and  $(-1, \xi , 0, 0 \ldots,  0)$, 
and  the centers of all
but one of the original unit spheres.   These centers form  the vertices 
of a regular $(n-1)$-simplex, so 
these $n-1$ ``orthogonal'' planes are equiangular
satisfying \eqn{N504a}, where $\theta$ is the angle between 
the normals of two facets of a regular $n$-simplex.
It follows that these orthogonal planes satisfy 
\eqn{N703}.  The
final ``orthogonal'' sphere through $T$ is orthogonal to the plane 
$x_1 = 1$ and all the $n$ original unit spheres.    Its center is 
thus $(1,0,0,\ldots)$ and its radius is $\xi$.  Hence it is also equiangular
with the $n-1$ ``orthogonal'' planes, with $\cos \theta_n =- \frac {1}{n-1}.$
These angles are all equal to the one formed by connecting the
vertices of a regular simplex to its center, i.e. the angle in
a triangle of sides $\xi, \xi$ and $2$.  
Finally, the last two ``orthogonal''  spheres meet at the
same angle in the plane $x_1 = 0.$
 ~~~$\bsq$
%
%

\section{Concluding remarks}
\setcounter{equation}{0}

This paper studied generalizations of the basic properties
of two-dimensional Apollonian packings and super-Apollonian
packings. We gave fairly complete answers but left open a
few problems.

The first problem  is to determine  a presentation of the
super-Apollonian group $\sA_n^{S}$ in dimensions $n \ge 3$, in terms of 
the given generators. Is this group always a hyperbolic Coxeter
group?

There are some unanswered questions concerning rational
and integral structures on Descartes configurations.  
In \S3.1 we raised the question of 
determining those dimensions $n$ in which 
the Descartes form $Q_{D,n}$, the Wilker form $Q_{W, n}$ 
and the Lorentz form $Q_{\sL, n}$ are rationally equivalent.
In \S6 we gave a necessary and
sufficient condition for rational equivalence of the
pair (Descartes, Wilker), namely that $n=2k^2$
or $n=(2k-1)^2$.  For the other 
two pairs (Descartes, Lorentz) and
(Wilker, Lorentz),  we gave  the necessary conditions
$n=2k^2$ and $n=2k$, respectively. It remains to determine
necessary and sufficient conditions in these cases.  

Finally we left open the question of whether there 
is any dimension $n \ge 3$ in which there  exist
strongly integral Descartes configurations.
Although strong integrality will not be
preserved under the action of the super-Apollonian group
(since  $n \ge 4$), it would be 
preserved under the action of the dual Apollonian group 
$\sA_n^{\perp}$, 
which consists of integer matrices in all dimensions.

%
%
%
\newpage
\section{Appendix.  M\"{o}bius Group Action}
\setcounter{equation}{0}
The  M\"{o}bius group action given in
Appendix A of part I  straightforwardly extends to
$n$-dimensions. 

The {\em (general) M\"{o}bius group} M\"{o}b$(n)$ is the group generated by 
reflections in spheres or planes in the one-point compactification
$\hat{\RR}^n = \RR^n \cup \{ \infty\}$ of $\RR^n$.
see Beardon \cite[Chapter 3]{Be83}. 
(He denotes it $GM(\hat{\RR}_n)$.)
This group 
has two connected components, and we let
$\mbox{M\"{o}b}(n)_{+}$ denote the connected component of
the identity. 
The {\em extended M\"{o}bius group} $GM^{*}(n)$ is 
defined by $GM^{\ast}(n) := \mbox{M\"{o}b}(n) \times \{ -\bI, \bI\}.$
Here $\{ -\bI, \bI\}$ are in the center of this group, and we
write elements of $GM^{\ast}(n)$ as $\pm \fg$, in which 
$\fg \in \mbox{M\"{o}b}(n)$, and the sign indicates which of
$\pm \bI$ occurs. 
The group $GM^{\ast}(n)$ has four connected
components. 

The purpose of this Appendix is to  define 
an action of $GM^{*}(n)$ on the right on the parameter  space $\sM_{\dd}^n$,
given in Theorem~\ref{thA1} below. This amounts to
finding an explicit isomorphism between $GM^{*}(n)$
and $Aut(Q_{W, n})$. The case $n=2$ was treated in
Appendix A of part I (\cite[Theorem 7.2]{GLMWY21}), was shown that
the M\"{o}b$(2)$ action preserves (total) orientation of Descartes
configurations.  The same property holds for  
the M\"{o}b$(n)$ action treated here, by a similar proof
which we omit.

Relevant isomorphisms are pictured in Table \ref{tableA1} below.
The isomorphism between $GM^{*}(n)$
and $Aut(Q_{W, n})$
appears as the  horizontal
arrow  on the left in the top row.
This map when restricted to the smaller groups
M\"{o}b(n) and M\"{o}b$(n)_{+}$ give the other
two horizontal isomorphisms on the left side of
the table. 
Table \ref{tableA1}   also indicates isomorphisms on its right side  
to the
orthogonal group $O(n+1,1)$ and 
corresponding  subgroups, which we defer discussing 
until  after the following result.

\begin{theorem}\label{thA1}
Let $GM^{*}(n) := \mbox{M\"{o}b}(n) \times \{ \bI, - \bI\}$.
There is a unique 
isomorphism $\pi: GM^{*}(n) \to Aut(Q_{W, n})$,
with image elements $ \bV_{\pm \fg}:=\pi(\pm \fg)$,
such that the following hold.
 
(i) For  $\fg \in  \mbox{ M\"{o}b}(n)$
the augmented curvature-center coordinates
for each ordered, oriented Descartes configuration $\sD$ 
satisfy 
\beql{N9.2}
\bW_{\fg(\sD)} = \bW_{\sD} \bV_{\fg}^{-1}.
\eeq

(ii) The action of $-\bI$ on augmented
curvature-center coordinates is
\beql{N9.2b}
\bW_{-\sD} = \bW_{\sD}\bV_{-\bI}^{-1}= -\bW_{\sD}.
\eeq 
This reverses the orientation of the Descartes configuration.
\end{theorem}

\paragraph{Proof.}
We compute the action of 
$\mbox{M\"{o}b}(n)$
on augmented curvature-center coordinates.
Let $(\bar{b}, b, w_1, w_2, \dots, w_n)= (\frac{(\sum_{i=1}^n x_i^2) -r^2}{r}, 
\frac{1}{r}, \frac{x_1}{r}, \frac{x_2}{r}, \dots, \frac{x_n}{r} )$
be the augmented curvature-center coordinates of the sphere
$$ \sum_{i=1}^n (y_i-x_i)^2 = r^2. $$
This sphere can be recovered from these coordinates via
\beql{281}
\sum_{i=1}^n (by_i - w_i)^2  = 1,
\eeq
and the orientation of the sphere (inside versus outside) is
determined by the sign of $b$. An oriented ``sphere at infinity''
is a hyperplane given by 
\beql{282}
 \by \cdot \bh = m,
\eeq
and its associated curvature-center coordinates are 
\beql{283}
(\bar{b}, b, w_1, w_2,\dots, w_n)= (2m, 0, h_1, h_2, \dots, h_n),
\eeq
where $\bh=(h_1, h_2, \dots, h_n)$ is the unit normal vector, and  the
orientation is given by the convention that the  
normal $\bh$ points inward. 

The group M\"{o}b(n) is generated by 

 (1) translations $\ft_{\by_0}(\by) = \by + \by_0$;

 (2) dilations $\fd_r(\by) = r \by$ with $r \in \RR$, 
$r > 0$;  

 (3) the rotation $\fo (\by) = \bO\by$, where $\bO$ is an $ n\times  n$ 
orthogonal matrix; 

 (4)  the inversion in the unit circle
$\fj_C(\by) =\frac{\by}{|\by|^2}.$

Given $\fg \in \mbox{M\"{o}b}(n)$,
we will let $\tilde{\fg}$ denote the 
corresponding action on  the curvature-center coordinates of an oriented
circle. 
The action of translation by $\by_0=(y_{0,1}, y_{0,2}, \dots, y_{0,n})$ is
\begin{multline} \label{284} 
\tilde{\ft}_{\by_0}(\bar{b}, b,  w_1  , w_2, \dots, w_n)= \nonumber \\ 
(\bar{b}+2\sum_{i=1}^n w_iy_{0,i} +b\sum_{i=1}^n y_{0,i}^2, b, w_1 + by_{0,1},
 w_2 + by_{0,2}, \dots, w_n + by_{0,n}).
\end{multline}
The action of a dilation with $r \in \RR$, $(r >0)$ is given by 
\beqs
\tilde{\fd}_\lambda (\bar{b}, b, w_1, w_2, \dots, w_n) = 
(r\bar{b}, b/r,~ w_1,  w_2, \dots, w_n). 
\eeqs
The action  of rotation $\fo$ with orthogonal matrix $\bO$ is
\beqs
\tilde{\fo}(\bar{b}, b, w_1, w_2, \dots, w_n) = (\bar{b}, b, w_1',  w_2', \dots, w_n'),
\eeqs
where $(w_1', w_2', \dots, w_n') =(w_1, w_2, \dots, w_n) \bO^T$. 
The action of inversion in the unit circle is
$$
\tilde{\fj}_C( \bar{b}, b, w_1, w_2, \dots, w_n) = (b, \bar{b}, w_1, w_2,\dots, 
w_n).
$$
All of these actions apply
to ``spheres at infinity'' and extend  to linear maps on the $(n+2) 
\times (n+2)$ matrices $\bW_{\sD}$. 

The translation operation is given by right multiplication by
the matrix 
\beqs
\bV_{\ft_{\by_0}}^{-1} := \left[
\begin{array}{ccccc}
1 &  0  &  0  & \ldots & 0 \\
\sum_{i=1}^n y_{0,i}^2  &  1   &  y_{0,1}   & \ldots &  y_{0,n}  \\
2y_{0,1}  &  0   & \ddots    &  &  \\
\vdots    &  \vdots &   & \bI_n &   \\
2y_{0,n} & 0 &   &  &  \ddots
\end{array}
\right],
\eeqs
and one verifies  \eqref{N9.2} holds by direct computation.
For the  dilation $\fd_r$, with $r \in \RR$,
$r > 0$  the right action is by the matrix
\beqs
\bV_{\fd_r}^{-1} :=\left[
\begin{array}{ccc}
r         &   0            &     0     \\
0         & 1/r            &     0     \\
0         & 0              &   \bI_n  
\end{array} 
\right].
\eeqs
For rotation $\fo$, the  right action is by the matrix 
\begin{equation*}
\bV_{\fo}^{-1} = 
:=\left[
\begin{array}{ccc}
1         &  0             &     0    \\
0         &  1             &     0     \\
0         &  0             & {\bf O}^T  
\end{array} 
\right].
\end{equation*}
For  the inversion $\fj_C$ in the unit circle, the permutation matrix 
\begin{equation*}
\bV_{\fj_C}^{-1} =\bV_{\fj_C}= \bP_{(12)} = \left[
\begin{array}{ccc}
0         &   1            &     0     \\
1         &   0            &     0     \\
0         &   0            &   \bI_n 
\end{array} 
\right].
\end{equation*}

It is easy to verify that the above matrices are all in $Aut(Q_W)$.
Since $\bQ_{W,n}=\bA^T \bQ_{\sL,n} \bA $ where 
$$
\bA=\sqrt{2} \left[
\begin{array}{ccc}
1  & 1 & 0  \\
-1 & 1 & 0 \\
0  & 0 &  \bI_n 
\end{array} 
\right],
$$
we have that a matrix $\bU \in Aut(Q_{W,n})^\uparrow$ if and only if
$\be_1^T \bA \bU \bA^{-1} \be_1 >0$. One checks that all the 
above matrices are actually in $Aut(Q_{W,n})^\uparrow$, 
so that the map  so far defines  a homomorphism of $\mbox{M\"{o}b}(n)$
into $Aut(Q_{W,n})^{\uparrow}\simeq O(n+1, 1)^{\uparrow}$,
identified with the  isochronous Lorentz group.
The group $\mbox{M\"{o}b}(n)$ acts simply transitively on
ordered Descartes configurations, as observed by Wilker
\cite[Theorem 3, p. 394]{Wi81},
and the group $Aut(Q_{W,n})$ acts simply transitively on
ordered, oriented Descartes configurations, as implied by Theorem 3.1.
Because  $Aut(Q_{W,n})^{\uparrow}$ is of index $2$ in 
$Aut(Q_{W,n}) \simeq O(n+1,1)$,
we conclude that the map so far defines  an isomorphism 
of $\mbox{M\"{o}b}(n)$  onto $Aut(Q_{W,n})^{\uparrow}$.

To complete the proof, we define  the action of $-I$ to be 
\begin{equation}~\label{flip}
(\bV_{-I})^{-1} = \bV_{-I} = - \bI_{n+2} .
\end{equation}
It has the effect of reversing  (total) orientation of the 
Descartes configuration, and does not correspond to 
a conformal transformation. Since $-\bI_{n+2}\notin Aut(Q_{W,n})^{\uparrow}$,
adding it gives the desired isomorphism of
$GM^*(n)$ onto $Aut(Q_{W,n})$.
~~~$\bsq$ \\

%
%

\begin{table}[t]
 \begin{tabular}[t]{|c|c|}
 \hline
  \begin{diagram}  
   GM^{*}(n)=  M\ddot{o}b(n)\times \{-\bI, \bI\}    & \rTo^\sim & Aut(Q_{W, n})  & 
         \rTo^\sim  &   {O(n+1,1)}  \\
   \dTo^\pi \uInto &   &  \dTo^\pi \uInto &    & \dTo^\pi \uInto  \\
   M\ddot{o}b(n)  & \rTo^\sim & 
Aut(Q_{W, n})^\uparrow 
         &  \rTo^\sim & O(n+1,1)^\uparrow  \\ 
   \uInto &  &  \uInto& & \uInto  \\
   M\ddot{o}b(n)_{+}    & \rTo^{\sim} 
         & Aut(Q_{W, n})^\uparrow_+ &  \rTo^\sim  &  O(n+1,1)^\uparrow_+  \\ 
  \end{diagram}
 &  
  \begin{diagram} 
  \text{Orthogonal group} \\ 
   {} \\
  {\text{Orthochronous} \atop  \text{orthogonal group}} \\ 
   {}\\
  {\text{Special orthochronous}  \atop \text{orthogonal  group}} \\
 \end{diagram} 
 \\
 \hline
\end{tabular}
\caption{Group Isomorphisms}~\label{tableA1}
\end{table}

We now return to the data in  Table ~\ref{tableA1}, 
giving  the isomorphisms of $Aut(Q_{W, n})$ and its
subgroups to the orthogonal group $O(n+1,1)$ and its
two subgroups 
$O(n+1,1)^{\uparrow}$ the
orthochronous orthogonal  group, and $O(n+1,1)_{+}^{\uparrow}$
the special  orthochronous orthogonal group,
which is the connected component of the identity
of the orthogonal group  $O(n+1,1)$.
The set of isomorphisms
given by the three horizontal arrows
on the right in Table~\ref{tableA1}
are obtained by any fixed  choice of real matrix $\bA$
that intertwines $Q_{W, n}$ and $Q_{\sL, n}$ by
$\bQ_{W} = \bA^T \bQ_{\sL, n} \bA$, in which case the 
isomorphism is $Aut(Q_{W, n}) = \bA^{-1} O(n+1,1) \bA$
sending $\bV \mapsto \bA \bV \bA^{-1}$.
Such matrices $\bA$  exist in all dimensions. 
It is shown in \S3.1 that rational matrices $\bA$ do
not exist in all dimensions; a necessary condition
for existence is that the dimension $n$ be even.

{\tt
\begin{tabular}{lllll}
email: & graham@ucsd.edu \\
& jagarias@umich.edu \\
  & colinm@research.avayalabs.com \\
 & allan@research.att.com \\
& cyan@math.tamu.edu
\end{tabular}
 }

\end{document}